\documentclass[leqno,12pt]{article}%
\usepackage{amsmath,amssymb}
\usepackage{amsthm}
\usepackage{dsfont}
\usepackage{amsfonts}%
\usepackage{color}

\newtheorem{Theorem}{\hspace{\parindent}\bf Theorem}[section]
\newtheorem{Lemma}{\hspace{\parindent}\bf Lemma}[section]

\newtheorem{Corollary}{\hspace{\parindent}\bf Corollary}[section]

\DeclareMathOperator*{\essinf}{ess\,inf}
\DeclareMathOperator*{\esssup}{ess\,sup}

\begin{document}

\title{\textbf{Singular integral equations with applications to travelling waves for doubly nonlinear diffusion}}
\author{Alejandro G\'arriz}

\maketitle

\begin{abstract}
We consider a family of  singular Volterra integral equations that appear  in the study of monotone travelling-wave solutions for a family of diffusion-convection-reaction equations involving the $p$-Laplacian operator. Our results extend the ones due to B.\,Gilding for the case $p=2$. The fact that $p\neq2$ modifies the nature of the singularity in the integral equation, and introduces the need to develop some new tools for the analysis. The results for the integral equation are then used to study the existence and properties of travelling-wave solutions for doubly nonlinear diffusion-reaction equations in terms of the constitutive functions of the problem.
\end{abstract}


\vskip 5cm

\noindent{\makebox[1in]\hrulefill}\newline2010 \textit{Mathematics Subject
Classification.}
45G05,  
45D05,  
35C07,  
35K57,  
35K59.   

\newpage
\section{Introduction}
\label{sect-introduction} \setcounter{equation}{0}

The goal of this paper is to study the Volterra integral equation
\begin{equation}\label{eq:formal_eq}
 x(t) =f(t) + \int_0^t \frac{g(s)}{x^\alpha (s)}\, {\rm d}s,\quad x(t)\ge 0,\quad t\geq0,
\end{equation}
where $\alpha\in \mathbb{R}_+:=(0,\infty)$, and
\begin{equation*}
\label{eq:conditions.alpha.f.g}
f\in C(\overline{\mathbb{R}_+}),\quad f(0)\geq 0, \quad
g\in L^1((0,\tau))\text{ for all }\tau\in\mathbb{R}_+,
\end{equation*}
extending the results obtained by Gilding in \cite{Gilding-1993} for the case $\alpha=1$.  We are mainly interested in the singular case $f(0)=0$. As shown in~\cite{Gilding-Kersner-PRSE-1996}, this equation arises in the study
of travelling waves in nonlinear reaction-convection-diffusion
processes
\begin{equation}\label{eq:complete_equation}
u_t=\Delta_p (a(u)) +  (b(u))_x + c(u),\quad x\in\mathbb{R},\, t\geq 0,
\end{equation}
involving the $p$-Laplacian operator,
$$
\Delta_p u=\mathop{\rm div}(|\nabla u|^{p-2}\nabla u), \quad p>1.
$$
\emph{Travelling waves or fronts} are special solutions of the form $u_\sigma(x,t)=V_\sigma(x-\sigma t)$ for some speed $\sigma$ and profile, that may depend on the speed, $V_\sigma$. They  are important at least in relation with two aspects: finite speed of propagation and large time behaviour.

Equation~\eqref{eq:complete_equation} is said to display \emph{finite speed of propagation} if a non-negative solution which has bounded support with respect to the spatial variable at some initial time also presents this property at later times. As proved in~\cite{Gilding-Kersner-PRSE-1996} (see also~\cite{Gilding-Kersner-JDE-1996}), equation~\eqref{eq:complete_equation} displays finite speed of propagation if and only if it possesses a travelling-wave solution whose profile's support  is bounded from above.  For this purpose it is enough to have a \emph{local} travelling wave, with a profile only defined in a half-interval~$(\omega,\infty)$.

A bounded \emph{global} (defined in the whole real line) travelling-wave profile $V_\sigma$ that is monotonic, but not constant, and such that
$$
V_\sigma(\xi)\to \ell^-\text{ as } \xi\to -\infty\quad\text{and } V_\sigma(\xi)\to \ell^+\text{ as } \xi\to \infty
$$
for some $\ell^-,\ell^+\in\mathbb{R}$ such that $c(\ell^-)=c(\ell^+)=0$ is said to be a \textit{wavefront profile} from $\ell^-$ to $\ell^+$ (the associated travelling wave is a \emph{wavefront solution}).  Wavefronts are already known to be important in the description of the large-time behaviour of more general solutions for wide classes of initial data when $p=2$, $a(u)=u^m$, $m\ge1$, and $b\equiv0$ for several reaction nonlinearities $c$; see for instance~\cite{Fisher-1937, Kolmogorov-Petrovsky-Piscounov-1937, Fife-McLeod-1975, Aronson-Weinberger-1978, Uchiyama-1978,   Bramson-1983, Biro-2002,  Du-Quiros-Zhou-Preprint,  Garriz-preprint}. They are also expected to give the large time behaviour when $p\neq2$. As a first step in this direction we have the papers~\cite{Audrito-Vazquez-2017-JDE, Audrito-2019}.

Let us explain the connection between the integral equation~\eqref{eq:formal_eq} and travelling waves for~\eqref{eq:complete_equation} in more detail, a connection that was first explored for the case $p=2$ in~\cite{Gilding-1996}.  Let $u_\sigma(x,t) = V_\sigma(\xi)$, $\xi =x-\sigma t$, be a travelling-wave solution to~\eqref{eq:complete_equation}  taking values in $\bar I$, where $I=(0,\ell)$, $0<\ell\le \infty$.  The profile $V$ (we drop the subscript for simplicity) satisfies the ordinary differential equation
\begin{equation}\label{eq:differential.equation}
\Delta_p (a(V)) + (b(V))_\xi + c(V) + \sigma V_\xi=0
\end{equation}
in a weak sense; see Section~\ref{sect-tw.general.results} for a precise definition.
Taking $g:= -(a(V))_\xi = -a^\prime(V)V_\xi$, we arrive to the system of equations
$$
V_\xi = \frac{-g}{a^\prime (V)}, \qquad
            g_\xi = \frac{1}{(p-1)|g|^{p-2}} \left( c - \frac{g(b^\prime + \sigma)}{a^\prime}\right).
$$
The function $a$ is assumed to be nondecreasing, so that  equation~\eqref{eq:complete_equation} is parabolic (may be degenerate or singular). Hence, if we restrict to nonincreasing profiles, we may assume that $|g| = g$. Thus, the trajectories in the phase plane of this system satisfy
$$
\frac{{\rm d}g}{{\rm d}V} = \frac{1}{p-1}\left( \frac{b^\prime + \sigma}{g^{p-2}} - \frac{a^\prime c}{g^{p-1}}\right).
$$
If~\eqref{eq:differential.equation} admits a solution $V$ such that $V \to 0$, and $(a(V))_\xi \to 0$ as $\xi \to \infty$, this $V$ is represented by a trajectory that approaches the point $(V,g)=(0,0)$. Thus, setting
$$
\theta(s) = g^{p-1}(\xi),\qquad s=V(\xi),
$$
and integrating the equation of the trajectories through (0,0) we arrive to
\begin{equation}\label{eq:complete_integral_equation}
\theta(s) = b(s)+\sigma s - \int_0^s \frac{a^\prime (r) c(r)}{\theta^{\alpha} (r)}\, {\rm d}r,\qquad \alpha = 1/(p-1)>0,
\end{equation}
which lies within the class of integral equations of the form~\eqref{eq:formal_eq}. Note that in terms of the original problem $\theta$ represents the flux as a function of the height of the solution.


Conversely, if the flux $\theta$  satisfies the \emph{integrability condition}
\begin{equation}\label{eq:integrability_condition}
\int_{s_0}^{s_1} \frac{a^\prime (r)}{\theta^\alpha(r)}\, {\rm d}r <\infty \quad \text{for all } 0<s_0<s_1<\delta,
\end{equation}
for $\delta=\ell$, a solution of~\eqref{eq:differential.equation} can be constructed from a solution of~\eqref{eq:complete_integral_equation} by means of
\begin{equation}
\label{eq:construction}
\int_\nu^{V(\xi)} \frac{a^\prime (r)}{\theta^\alpha(r)}\, {\rm d}r = \xi_0 - \xi
\end{equation}
for some $\nu$ in the domain of definition of $\theta$ and some $\xi_0\in \mathbb{R}$. Let us remark that the profile obtained in this way is strictly monotone in $\{\xi\in\mathbb{R}: 0<V(\xi)<\ell\}$.

%

It turns out that the existence of wavefront profiles from $\ell$ to $0$ is equivalent, if~$\ell$ is finite, to the existence of global solutions to~\eqref{eq:complete_integral_equation} satisfying the integrability condition.
\begin{Theorem}\label{thm:correspondencia_ecuaciones_integral_y_diferencial}
Let $\ell<\infty$. The differential equation~\eqref{eq:differential.equation} has a wavefront profile from $\ell$ to $0$ with speed $\sigma$ if and only if the integral equation~\eqref{eq:complete_integral_equation} with $\alpha=1/(p-1)$ has a solution $\theta$  such that $\theta(\ell) = 0$ satisfying the integrability condition~\eqref{eq:integrability_condition} on $\bar{I}$.
\end{Theorem}
We omit the proof, since it is similar to the one for the case $p=2$ (that is, $\alpha=1$) given in~\cite{Gilding-Kersner-2004}.

Let us remark that there may be several wavefront profiles  leading  to the same solution $\theta$ of the integral equation. For instance, equation
$$
\displaystyle u_t=\Delta_pu +
\begin{cases}
(p-1)(2u - 4u^2)^{\frac{p-2}{2}}2^p u^p (8u-3)  &\text{for }0\leq u\leq1/2,\\[8pt]
(p-1)(2(1-u) - 4(1-u)^2)^{\frac{p-2}{2}}2^p(1-u)^p (8u-5) &\text{for }1/2< u\leq 1,
\end{cases}
$$
with $p>2$ admits a family of wavefront profiles from 1 to 0 with speed $\sigma=0$ given by
$$
V(\xi)=\begin{cases}
\frac{1+2(\xi_0-\xi)^2}{2(1+(\xi_0-\xi)^2)} &\text{for }\xi<\xi_0, \\[6pt]
1/2  &\text{for }\xi_0\leq \xi\leq \xi_1, \\[6pt]
\frac1{2(1+(\xi-\xi_1)^2)}  &\text{for }\xi_1<\xi,
\end{cases}
$$
for any pair or real numbers $\xi_0\le\xi_1$.  The corresponding flux
$$
\theta(s)= \left\{ \begin{array}{lcc}
\left(2s\sqrt{2s(1-2s)}\right)^{p-1} &\text{for }0\leq s\leq 1/2, \\[8pt]
\left(2(1-s)\sqrt{2(1-s)(2s-1)}\right)^{p-1}  &\text{for }1/2<s\leq 1,
\end{array}
\right.
$$
satisfies the integral equation~\eqref{eq:complete_integral_equation} regardless the values of $\xi_0$ and $\xi_1$.

On the other hand, starting from this solution $\theta$ of the integral equation and using~\eqref{eq:construction} with $\nu\in(0,1)$ we obtain the wavefront profiles
\begin{equation}
\label{eq:profiles.example}
V(\xi)=\begin{cases}
\frac{1+2(\eta-\xi)^2}{2(1+(\eta-\xi)^2)} &\text{for }\xi\leq\eta, \\[6pt]
\frac1{2(1+(\xi-\eta)^2)}  &\text{for }\eta<\xi,
\end{cases}
\qquad\text{where }
\eta=\begin{cases}
-\left(\frac{1-2\nu}{2\nu}\right)^{1/2}&\text{for }\nu\leq 1/2, \\[8pt]
\left(\frac{2\nu -1}{2(1-\nu)}\right)^{1/2}  &\text{for }\nu>1/2.
\end{cases}
\end{equation}
We observe that in this example the flux vanishes only at heights where the reaction vanishes, $s=0,1/2,1$. This is true not only in this particular case, but in general, as can be shown by adapting the proof given in~\cite[Lemma 2.24]{Gilding-Kersner-2004} when $\alpha=1$. Note, however, that the opposite implication does not hold: the reaction term may vanish at points where the flux does not.

When $p=2$ the  example above coincides with the one given in~\cite[Application 2.22]{Gilding-Kersner-2004}. Nevertheless, in sharp contrast with that case, if $p>2$ the reaction term is continuous in $(0,\ell)$. Moreover,  if $p\ge 4$, $c$ and hence $ca'$ are differentiable in the same open interval. However, the flux~$\theta$ vanishes at a point in $(0,\ell)$, namely $s=1/2$. This is not possible when $\alpha=1$, since for that value of the parameter the differentiability of $ca'$ implies the positivity of the flux; see~\cite[Lemma~2.40]{Gilding-Kersner-2004}. The possibility of having a vanishing flux at positive heights for smooth reactions is introduced by the degeneracy of the $p$-Laplacian for $p>2$. In the singular case $1<p<2$ this possibility does not exist, as can be easily proved in the same way as for $\alpha=1$.

We also observe that all the wavefront profiles~\eqref{eq:profiles.example} arising from the integral equation are translates one of each other, with $\nu$ giving the point where~$V$ takes the value 1/2. These wavefront profiles correspond to the  ones with $\xi_0=\xi_1=\eta$, which are strictly monotonic.

Two wavefront profiles  are said to be \textit{indistinct} if one is a translation of the other. Otherwise they are \textit{distinct}. It turns out that there is a one-to-one correspondence between solutions of the integral equation~\eqref{eq:complete_integral_equation} in $[0,\ell]$ with  $\ell$ finite and the set of distinct wavefront profiles from $\ell$ to~$0$ of~\eqref{eq:complete_equation} that are strictly monotone in their positivity set.

\begin{Theorem}\label{thm:es_el_2.25_del_libro}
Let $\ell<\infty$. To every solution $\theta$ of~\eqref{eq:complete_integral_equation} satisfying~\eqref{eq:integrability_condition} in $[0,\ell]$ and $\theta(\ell)=0$ there corresponds precisely one distinct wavefront profile $V$ of equation~\eqref{eq:complete_equation} with speed $\sigma$ from $\ell$ to $0$ which  is strictly monotonic where it is positive. Moreover, there corresponds no other (non-strictly monotonic) distinct wavefront profile with speed $\sigma$ from $\ell$ to 0 if and only if $\theta>0$ in $(0,\ell)$.
\end{Theorem}
We omit the proof, since it is similar to the one of \cite[Theorem 2.26]{Gilding-Kersner-2004}.

\medskip

Back to the integral equation~\eqref{eq:formal_eq}, Gilding and Kersner have already considered it for $\alpha\neq1$ in~\cite{Gilding-Kersner-PRSE-1996} for the particular case 
$$
f(t)=\sigma t+\mu t^\gamma,\quad  \mu\in\mathbb{R}, \ \gamma>0,\qquad g(t)=c_0 t^\beta,\quad c_0\in\mathbb{R},\ \beta>-1. 
$$
They characterize there for which values of the parameters $\alpha,\mu,\gamma,c_0,\beta$ there is a local solution of the integral equation for some large enough value $\sigma$. This can be translated into the existence of a local travelling wave with a profile whose support is bounded from above, which is enough for their purposes.

In the present paper we obtain both necessary conditions and sufficient conditions in terms of $f$, $g$ and $\alpha$ for the existence of local solutions to~\eqref{eq:formal_eq}, in the spirit of~\cite{Gilding-1993}, which are enough to analyze the question of existence of travelling waves and their behavior for diffusion-reaction equations for reaction nonlinearities with none or one sign change. Though we follow ideas from~\cite{Gilding-1993}, the new degeneracy/singularity represented by the power $\alpha$ introduces the need to develop nontrivial generalizations of several tools that in addition offer more insight on the analysis.

Once the study of the integral equation is finished, we proceed to apply it to the study of travelling-wave solutions of~\eqref{eq:complete_equation}, following the ideas developed by Gilding and Kersner in the remarkable book~\cite{Gilding-Kersner-2004} for the case $p=2$. Such analysis when  $p\neq2$ has already been considered in~\cite{Enguica-Gavioli-Sanchez-2013,Gavioli-Sanchez-2015} for $a(u)=u$, and in~\cite{Audrito-Vazquez-2017-JDE, Audrito-2019} for $a(u)=u^m$ using a different technique. However, these papers have several restrictions on the reaction nonlinearity $c$.

\noindent\textsc{Organization of the paper. } Section~\ref{sect-integral_equation} represents the core of the article, where we study the integral equation~\eqref{eq:formal_eq} in detail. In Section~\ref{sect-tw.general.results} we exploit the connection between this integral equation and the diffusion-convection-reaction equation~\eqref{eq:complete_equation} to obtain some general results for the latter. We start by studying a local version of wavefronts, called semi-wavefronts, and then proceed to extend them to the whole real line. Finally, in Section~\ref{sect-travelling_waves_reaction_diffusion} we make full profit of the results of Section~\ref{sect-integral_equation} by applying them to study travelling waves of diffusion-reaction equations (no convection).

\section{The integral equation}
\label{sect-integral_equation} \setcounter{equation}{0}

Throughout this section we will study equation~\eqref{eq:formal_eq} in detail in the same fashion as Gilding did in his work \cite{Gilding-1993} when $\alpha=1$.

Let us define
$$
I(s,x)=\begin{cases}
            g(s)/x^\alpha \quad &\text{if}\ x>0, \\
            -\infty \quad &\text{if}\ g(s)<0\ \text{and}\ x=0, \\
            0 \quad &\text{if}\ g(s)=0\ \text{and}\ x=0, \\
            \infty \quad&\text{if}\ g(s)>0\ \text{and}\ x=0.
             \end{cases}
$$
Considering the integral in~\eqref{eq:formal_eq} as an improper Lebesgue integral, we will say that a function $x$ is a solution of equation~\eqref{eq:formal_eq} on $[0,\tau)$ if it is defined, real, nonnegative and continuous in $[0,\tau)$, $G(s,x(s))\in L^1_{\rm loc}((0,\tau))$, and
$$
\int_0^t I(s,x(s))\, {\rm d}s := \lim\limits_{\varepsilon \downarrow 0} \int_\varepsilon^t I(s,x(s))\, {\rm d}s
$$
exists and satisfies
$$
x(t)=f(t) + \int_0^t I(s,x(s))\, {\rm d}s \quad\text{for all }t\in(0,\tau).
$$
One can check that a solution in this sense satisfies
$
x(0)=f(0)$.

We devote a first subsection to give some continuation, uniqueness, and existence results for~\eqref{eq:formal_eq} which are valid in general. We next pay attention in three different subsections to some particular instances of the integral equation that play an important role in the analysis of wavefront profiles, taking advantage of the sign of the kernel $g/x^\alpha$,  as Gilding did in~\cite[Sections 7, 8 and 9]{Gilding-1993}. It is here where the analysis becomes significantly different to Gilding's work, requiring the development of new ideas.

\subsection{General results}

Most of the results in this section will be presented without proofs for the sake of simplicity, since they are very similar to the ones in~\cite{Gilding-1993}. This similarity is supported in the fact that the function $x\rightarrow x^\alpha$ is monotonic and differentiable if $x>0$.

We start with a prolongability result.

\begin{Theorem}\label{Thm:extendibility}
Let $x$ be a solution of~\eqref{eq:formal_eq} on  $[0,T)$, $T<\infty$. Then the limits
$$
x(T):=\lim\limits_{t\to T} x(t),\qquad\int_0^T \frac{g(s)}{x^\alpha(s)}\, {\rm d}s:= \lim\limits_{t\to T} \int_0^t \frac{g(s)}{x^\alpha(s)}\, {\rm d}s
$$
exist and are finite and
$$
x(T)=f(T) + \int_0^T \frac{g(s)}{x^\alpha(s)}\, {\rm d}s.
$$
Moreover, either $x(T)=0$ or $x$ is continuously extendible as a solution of equation~\eqref{eq:formal_eq} onto a bigger finite interval $[0,T')$.
\end{Theorem}

There is also a comparison principle.

\begin{Theorem}\label{lemma:formal_comparison}
Let $x_1$ denote a solution of the equation
\begin{equation}
\label{eq:thm.comparison}
x_1(t)=f_1(t) + \int_0^t \frac{g_1(s)}{x_1^\alpha (s)}\, {\rm d}s
\end{equation}
on some finite interval $[0,\delta)$. Suppose that $f(0)>f_1(0)$, $(f-f_1)$ is nondecreasing on $[0,\delta)$ and $g(t)\geq g_1(t)$ a.e.~on $(0,\delta)$. If~\eqref{eq:formal_eq} has a solution $x$ on $[0,\delta)$ such that $x(t) >0$ for all $t\in [0,\delta)$, then $x(t) > x_1(t)$ for all $t\in [0,\delta)$.
\end{Theorem}
\begin{proof}
Suppose that the lemma is false. Then there exists a value $t^*\in(0,\delta)$ such that $x(t^*)=x_1(t^*)>0$ and $x(t)>x_1(t)$ for all $t\in[0,t^*)$.

Let $t_1\in[0,t^*)$ be such that $x^\alpha (s) - x_1^\alpha(s)< K(x(s)-x_1(s))$ for all $s\in[t_1,t^*)$. This constant $K$ must exist, since the function $x^\alpha$ is differentiable and hence Lipschitz in this interval. After this, we take $t_0\in [t_1,t^*)$ such that $x_1(t)>0$ for all $t\in[t_0,t^*]$ and
$$
\int_{t_0}^{t^*} \frac{|g(s)|}{x^\alpha (s)x_1^\alpha (s)}\, {\rm d}s \leq 1/2K.
$$

We finish as in \cite{Gilding-1993}: using the equations satisfied by $x$ and $x_1$ we see that
$$
x(t)-x_1(t) \leq \int_t^{t^*} g(s)\left(\frac{1}{x_1^\alpha(s)} - \frac{1}{x^\alpha(s)}\right)\, {\rm d}s \leq \frac{1}{2} \|x-x_1\|_{L^\infty (t_0, t^*)},\quad t\in[t_0,t^*],
$$
so that $\|x-x_1\|_{L^\infty (t_0, t^*)}=0$, contradicting that $x(t)>x_1(t)$ for all $t\in[0,t^*)$.
\end{proof}

We now pay attention to the number of solutions.

\begin{Theorem}\label{thm:formal_uniqueness}
Equation~\eqref{eq:formal_eq} has none, one or an uncountable number of solutions. Moreover, if $\essinf\limits_{0<t<\tau} g(t) \geq 0$  for some $\tau >0$, then it  has at most one solution in $[0,\tau)$.
\end{Theorem}

In the non-singular case $f(0)>0$ existence and uniqueness follow from standard theory for nonlinear Volterra integral equations; see for instance~\cite{Gripenberg-Londen-Staffans-1990}.

\begin{Theorem}\label{thm:formal_f(0)>0}
If $f(0)>0$ then equation~\eqref{eq:formal_eq} has a unique positive solution $x$ on an interval $[0,\tau)$ such that either $x(t)\to 0$ as $t \uparrow \tau$ or $\tau=\infty$.
\end{Theorem}

Hence,  we only have a difficulty when $f(0)=0$. The idea to deal with it is to lift the datum $f$ by a constant $\mu>0$ and then pass to the limit. Thus, we will construct a solution as $\lim_{\mu \downarrow 0}x(t;\mu)$, where $x(t;\mu)$ denotes the unique positive solution to
\begin{equation}\label{eq:formal_mu}
x(t)=\mu + f(t) + \int_0^t \frac{g(s)}{x^\alpha(s)}\, {\rm d}s,
\end{equation}
which has an interval of existence $[0,T(\mu))$.

By Lemma~\ref{lemma:formal_comparison} we have for all $0<\mu_1<\mu_2<\infty$ that $T(\mu_1)\leq T(\mu_2)$ and
$$
x(t;\mu_1)<x(t;\mu_2)\quad\text{for all } t\in[0,T(\mu_1)).
$$
Moreover, $T(\mu)\to\infty$ as $\mu\to\infty$. These properties allow us to define
$$
\tilde{x}(t;0)= \inf\{x(t;\mu):\mu\in \mathbb{R}_+ \text{ such that }T(\mu)>t\}.
$$
Following \cite[Section 4]{Gilding-1993} we can prove that whenever~\eqref{eq:formal_eq} has a solution, $\tilde x$ is a \emph{maximal} solution of the same equation.

\begin{Theorem}
If equation~\eqref{eq:formal_eq} has a solution $x$ on an interval $[0,\delta)$ then $\tilde{x}(t;0)$ is a solution of~\eqref{eq:formal_eq} on an interval $[0,\tau) \supseteq [0,\delta)$ and $\tilde{x}(t;0) \geq x(t)$ for all $t\in [0,\delta)$. Furthermore, either $\tau =\infty$ or $\tilde{x}(t;0)\to 0$ as $t\uparrow\tau$.
\end{Theorem}

We can now give  sufficient conditions for $\tilde x$ to be a solution, necessarily the maximal one, proceeding as in \cite[Section 5]{Gilding-1993}.
\begin{Theorem}\label{thm:formal_comparison_maximal}
Let $x_1$ be a solution of~\eqref{eq:thm.comparison}
on some finite interval $[0,\tau)$.
\begin{itemize}
\item[\rm (i)]  Suppose that $f(0)\geq f_1(0)$, $(f-f_1)$ is nondecreasing on $[0,\tau)$, $g(t)\geq g_1(t)$ a.e. on $(0,\tau)$ and $g/x_1^\alpha \in L_{\rm loc}^1((0,\tau))$. Then $\tilde{x}(t;0)$ solves~\eqref{eq:formal_eq} on $[0,\tau)$ and
$$
\tilde{x}(t;0)\geq x_1(t)\text{ for all }t\in [0,\tau).
$$

Moreover, if $(f-f_1)(s)<(f-f_1)(t)$ for all $s\in [0,t)$ for a fixed $t\in[0,\tau)$, then $\tilde{x}(t;0)=x_1(t)$ if and only if $\tilde{x}(t;0)=0$.

\item[\rm (ii)]  Suppose that $f(t)\geq f_1(t)$ for all $t\in[0,\tau)$ and $\min\{0,g(t)\}\geq g_1(t)$ a.e.\,on $(0,\tau)$. Then $\tilde{x}(t;0)$ solves~\eqref{eq:formal_eq} on $[0,\tau)$ and
$$
\tilde{x}(t;0)-f(t)\geq x_1(t)-f_1(t)\text{ for all }t\in [0,\tau).
$$
\end{itemize}
\end{Theorem}

\subsection{Existence: negative kernels}\label{sect-Negative_Kernel}

We now concentrate on the case of negative kernels,
\begin{equation}
\label{eq:assumption.negative.kernel}
\esssup\limits_{0<t<\tau} g(t)\leq 0\quad\text{for some }\tau\in(0,\infty].
\end{equation}
As before, we only need to consider the case $f(0)=0$.

Some special functions solving~\eqref{eq:formal_eq} for particular data $f$ will play an important role in the analysis. The first one,
\begin{equation}\label{eq:G}
G(t):=\left|(\alpha+1)\int_0^t g(s)\, {\rm d}s\right|^{\frac{1}{\alpha + 1}},
\end{equation}
solves~\eqref{eq:formal_eq} if $f\equiv0$.  We may  assume that
\begin{equation}
\label{eq:condition.positivity.G}
G(t)>0\quad \text{for all }t\in (0,\tau),
\end{equation}
since otherwise  our equation reduces to $x(t)=f(t)$ on a small neighbourhood of~0, a trivial case. The second kind of special solutions have the form $x(t)=zG(t)$, $z>0$, and solve the integral equation with $f(t)=kG(t)$  if and only if $z$ is a root of the \lq\lq fractional polynomial''
\begin{equation*}\label{eq:polinomio.+}
P^+_k(z):=z^{\alpha+1} - kz^\alpha + 1 =0, \qquad z>0.
\end{equation*}
This polynomial attains its minimum in $\mathbb{R}_+$ for $\displaystyle z=z_{min}:=\frac{k\alpha}{\alpha+1}$, and the number of its (positive) roots depends on the sign of $k-k_0$, where
$$
k_0:=(\alpha + 1) \alpha^{-\frac{\alpha}{\alpha +1}}.
$$
Indeed, a simple analysis shows that if $k>k_0$, there are two positive roots of~$P_k^+$, say $z_1$ and $z_2$, while if $k=k_0$ there is only one, $z_0$, and $z_0=z_{min}$. It is worth noting that in the case of positive kernels $g$, as considered in the next subsection, instead of the polynomial $P^+_k$ we have to deal with
\begin{equation*}
\label{eq:polinomio.-}
P^-_k(z):=z^{\alpha+1} - kz^\alpha -1 =0.
\end{equation*}
The change of sign of the zero-order term makes the analysis simpler in this case.

Let us finally define the important auxiliary quantities
\begin{equation*}
\label{eq:def.L.J}
L(t):=\frac1{|\ln G(t)|}, \quad J(t):=\frac1{|\ln L(t)|},\quad
\quad t\in[0,\tau).
\end{equation*}
Note that $|g|=G'G^\alpha$, $L'=G'L^2/G$ and $J'=L'J^2/L$ almost everywhere in $(0,\tau)$. The dependence on $t$ will be omitted in what follows when it is clear.

\begin{Theorem}\label{thm:existence_negative_kernel_1}
Assume that \eqref{eq:assumption.negative.kernel} and \eqref{eq:condition.positivity.G} hold. If
$$
f(t)\geq \left(k_0G - K_0GL^2(1+J^2)\right)(t)\quad\text{for all } t\in(0,\tau),\quad \text{with } K_0:=\frac{\alpha^\frac{1}{\alpha+1}}{2(\alpha+1)},
$$
then equation~\eqref{eq:formal_eq} has a maximal solution $\tilde{x}(t;0)$ on an interval $[0,\delta)\subset[0,\tau)$, and
$$
0\geq \tilde{x}(t;0) - f(t)\geq- \alpha^\frac{-\alpha}{\alpha+1}\big(G +\frac{\alpha}{\alpha+1}GL(1+J+J^2)\big)(t)  \quad\text{for all } t\in(0,\delta).
$$
\end{Theorem}
\begin{proof}
 Let $\delta>0$ be such that $G(t)<\exp(-1)$ for $t\in[0,\delta)\subset [0,\tau)$, and let us define  $a:=\alpha/(\alpha+1)$ and
$$
\displaystyle S:=G + aGL(1+J+J^2).
$$
A direct computation shows that
$$
D:=\frac{\partial S}{\partial G} = 1 + aL(1+J+J^2) + aL^2(1+J)(1+2J^2),
$$
and through the Taylor series of the function $(1+x)^n$ we also see that
$$
\begin{aligned}
D^{1/\alpha}=&1+\frac{a}{\alpha}(1+J+J^2)L + \frac{a}{\alpha}\left( (1+J)(1+2J^2) + \frac{1-\alpha}{2\alpha} a  (1+J+J^2)^2   \right)L^2 \\
&+ \mathcal{O}(L^3).
\end{aligned}
$$
It is easy then to check that the function
$$
x_1:=\frac{AG}{D^{1/\alpha}},\quad\text{where } A:=\alpha^{\frac{1}{\alpha+1}},
$$
is a solution of our equation provided that  $f_1:=x_1 + A^{-\alpha}S$. This is equivalent to
$$
\begin{aligned}
f_1=&k_0G-K_0GL^2(1+J^2)\\
&-D^{-1/\alpha}\left(k_0GD^{1/\alpha} - K_0GL^2(1+J^2)D^{1/\alpha} - AG -\frac{1}{A^\alpha}SD^{1/\alpha}  \right).
\end{aligned}
$$
The theorem follows easily once we prove that this last quantity between parentheses is positive. In order to do so, we group terms by its order.

The terms of order $G$ are multiplied by the factor $k_0-A-A^{-\alpha}$, but this amounts to 0 by definition. The terms of order $GL(1+J+J^2)$ are multiplied by the factor $k_0a\alpha^{-1} - aA^{-\alpha} - aA^{-\alpha}\alpha^{-1} $, but again this amounts to 0. Note that in order for these two quantities to be 0 we have not used the value of $a$.

Recall now that $(1+J+J^2)^2=1+2J+3J^2+2J^3+J^4$, and let us see the term of order $GL^2$. It is multiplied by
$$
a^2\left( \alpha^\frac{-1-2\alpha}{\alpha+1} - \frac{\alpha+1}{2} \right) + a \alpha^\frac{-\alpha}{\alpha+1} - K_0,
$$
and again this is 0, this time thanks to the definition of $a$. The terms of order $GL^2J$ and $GL^2J^2$ suffer the same fate. Finally, the term of order $GL^2J^3$ is multiplied by a factor of $2K_0>0$. The rest of the result now follows from Theorem~\ref{thm:formal_comparison_maximal}.
\end{proof}

\begin{Theorem}\label{thm:no_existence_negative_kernel}
Assume that \eqref{eq:assumption.negative.kernel} and \eqref{eq:condition.positivity.G} hold. Let $K_0$ be as in Theorem~\ref{thm:existence_negative_kernel_1}. If
$$
f(t)\leq (k_0G - \beta GL^2) (t) \quad\text{for all } t\in(0,\tau)
$$
for some $\displaystyle\beta>K_0$, equation~\eqref{eq:formal_eq} has no solution.
\end{Theorem}
\begin{proof}
We define
$$
H(t):=G^{-1}(t)L^{-1}(t),\quad Y(t):=-\int_0^t\frac{g(s)}{x^\alpha(s)}\, {\rm d}s\quad\text{for } t\in(0,\delta').
$$
Note that
$$
H'(t) = \frac{H(t)g(t)}{G^{\alpha +1}(t)}(L(t)+1)\quad\text{for almost every } t\in(0,\delta').
$$
Our goal is to obtain an absurd estimate for $Y$.  For a start, since $x\leq f\leq k_0G$, we obtain, substituting in the equation, that $x\leq (k_0-k_0^{-\alpha})G$. We can repeat this process again and again to obtain that $x\leq z_0G$, with
$$
z_0:=k_0-\frac{1}{\left(k_0-\frac{1}{\left(k_0-\frac{1}{\left(\cdots\right)^{\alpha}}\right)^{\alpha}}\right)^{\alpha}}.
$$
Note that $z_0$ satisfies $z_0=k_0-\frac{1}{z_0^\alpha}$,
which is precisely the root of~$P^+_k$, hence the notation. Therefore
\begin{equation}\label{eq:71}
Y(t)\geq -z_0^{-\alpha}\int_0^t g(s)/G^\alpha(s)\, {\rm d}s = z_0^{-\alpha}G(t)\quad\text{for all }t\in(0,\delta').
\end{equation}

On the other hand we have that
$$
Y(t)H(t)-Y(\varepsilon)H(\varepsilon) =\int_\varepsilon^t \big( Y'(s)H(s)+Y(s)H'(s)\big)\, {\rm d}s
$$
for any $0<\varepsilon<t<\delta$. But substituting $H, H',Y, Y'$ and using~\eqref{eq:71} we get
$$
Y(t)H(t) -z_0^{-\alpha}L^{-1}(\varepsilon) \geq \int_\varepsilon^t \left\{\frac{|g|}{L G}\left(\frac{1}{x^\alpha}+\frac{x(L+1)}{G^{\alpha+1}}\right) -\frac{f(L+1)|g|}{G^{\alpha + 2}L}\right\}\, {\rm d}s.
$$
Now the term between parentheses can be estimated by studying the function $x^{-\alpha} + Ax$ for $x\geq 0$ and a constant $A\geq 0$. This function attains at the point $x=(\alpha/A)^{1/(\alpha+1)}$ a minimum value of $k_0A^{\alpha/(\alpha+1)}$. Therefore, and using our hypothesis,
$$
Y(t)H(t) -z_0^{-\alpha}L^{-1}(\varepsilon) \geq \int_\varepsilon^t \left\{\frac{|g|k_0(L+1)^{\frac{\alpha}{\alpha+1}}}{L G^{\alpha+1}} -\frac{(k_0-\beta L^2)(L+1)|g|}{G^{\alpha + 1}L}\right\}\, {\rm d}s,
$$
and, after throwing away the appearing term $\beta L^3$, we arrive to
$$
Y(t)H(t) -z_0^{-\alpha}L^{-1}(\varepsilon) \geq \int_\varepsilon^t \frac{|g|}{L G^{\alpha+1}}\left(k_0((L+1)^{\frac{\alpha}{\alpha+1}} -(L+1)) + \beta L^2\right)\, ds.
$$
It is easy now to check, via Taylor series for example, that
$$
(L+1)^{\frac{\alpha}{\alpha+1}}\geq 1+\frac{\alpha}{\alpha +1}L - \frac{\alpha}{2(\alpha+1)^2}L^2.
$$
Therefore, recalling that $|g|=G'G^\alpha$ and that $z_0=k_0\alpha/(\alpha+1)$, we obtain
$$
Y(t)H(t) -z_0^{-\alpha}L^{-1}(\varepsilon) \geq \int_\varepsilon^t \frac{G'}{G}\left( \left(\beta-K_0\right)L + (k_0-z_0)  \right)\, {\rm d}s.
$$
Therefore, using that $L'=L^2G'/G$ and integrating,
$$
Y(t)H(t)  \geq   \left(\beta-K_0\right) \ln\left(\frac{L(t)}{L(\varepsilon)}\right)+ (k_0-z_0)\ln\left(\frac{G(t)}{G(\varepsilon)}\right) + z_0^{-\alpha}L^{-1}(\varepsilon).
$$
Since $k_0-z_0 = z_0^{-\alpha}$, taking $\varepsilon$ small enough so that $L^{-1}(\varepsilon) + |\ln(G(\varepsilon))|=0$, we get
$$
Y(t)H(t) \geq \left(\beta-K_0\right) \ln\left(\frac{L(t)}{L(\varepsilon)}\right)+ z_0^{-\alpha}\ln L^{-1}(t),
$$
and we get our contradiction by making $\varepsilon$ go to 0 whenever $\beta > K_0$.
\end{proof}

We can give an additional necessary condition for the existence of solutions of our integral equation. For it, we need the following lemma concerning non-negative kernels, which will be also used in the next subsection.

\begin{Lemma}\label{lemma:ode_for_integral_equation}
Suppose that $\essinf\limits_{t\in(0,\tau)} g(t)\geq 0$ and $f(t)\leq \beta G(t)$ for all $t\in [0,\tau)$ for some $\tau\in(0,\infty)$ and $\beta\in(-\infty,\infty]$. Then, given any $\mu>0$, there exists  $\rho>0$ such that
$$
x(t;\mu)\geq \rho+f(t)+\beta_0^{-\alpha}G(t)\quad\text{for all }t\in[0,\min\{T(\mu), \tau\}),
$$
where $\beta_0$ is the unique positive root of the \emph{fractional polynomial} $P^-_\beta$
when $\beta<\infty$ and $\beta_0^{-\alpha}:=0$ when $\beta=\infty$.
\end{Lemma}
\begin{proof}
The case $\beta=\infty$ follows easily from~\eqref{eq:formal_mu} simply by taking $\rho=\mu$.

For finite $\beta$ we argue by contradiction. Note that $x(0;\mu)=\mu>0$. Suppose that there exists $t_0 \in [0,\min\{T(\mu), \tau\})$ such that
\begin{equation}\label{eq:contrad}
x(t;\mu)>f(t) + \beta_0^{-\alpha}G(t) \quad \text{for all } t\in (0,t_0),\quad x(t_0;\mu)=f(t_0) + \beta_0^{-\alpha}G(t_0).
\end{equation}

We set
$$
Y(t):=\mu + \int_0^t \frac{g(s)}{x^\alpha(s;\mu)}\, {\rm d}s.
$$
If we multiply equation~\eqref{eq:formal_mu} by $Y'$ we get
\begin{equation}
\label{eq:edo}
0=YY' + fY' - x^{1-\alpha}G^\alpha G'
\end{equation}
for every $t\in(t_0,t_1)$. Substituting $x=f+Y$ and recalling \eqref{eq:contrad} and that $f\leq \beta G$,
$$
0\leq (\beta G + Y)^\alpha Y' - G^\alpha G'\text{ and }  Y(t)> \beta_0^{-\alpha}G(t)\text{ for all }t\in (0,t_0),
$$
with
$$
G(0)=0, Y(0)=\mu\text{ and }Y(t_0)= \beta_0^{-\alpha}G(t_0).
$$
From here it is easy to arrive to
$$
\frac{G(t_0)^{\alpha+1}}{\alpha+1} \leq \int_0^{t_0}(\beta G + Y)^\alpha Y'\,\rm{d}t < \int_0^{t_0}(\beta \beta_0^\alpha + 1)^\alpha Y^\alpha Y'\,\rm{d}t,
$$
and, recalling that $\beta_0^{\alpha+1} = \beta\beta_0^\alpha + 1$, we obtain
$$
\mu^{\alpha+1} < Y(t_0)^{\alpha+1} - \big(\beta_0^{-\alpha}G(t_0)\big)^{\alpha + 1} = 0,
$$
which is a contradiction with the fact that $\mu>0$.
\end{proof}

This lemma corresponds to~\cite[Lemma 11]{Gilding-1993} for the case $\alpha=1$, and yields the following theorem, arguing as in  Theorem~13 in the same paper.

\begin{Theorem} Let  $\beta$ and $\beta_0$  be as in the previous lemma, and assume~\eqref{eq:assumption.negative.kernel}. If
$$
\begin{array}{c}
f(t)<\beta_0^{-\alpha}G(t)$ for some $t\in(0,\tau)\text{ and}\\[6pt]
f(s)-f(t)\leq  \beta \{G^{\alpha +1}(t) - G^{\alpha +1}(s)\}^\frac{1}{\alpha + 1}\quad\text{for all }0\leq s<t,
\end{array}
$$
then equation~\eqref{eq:formal_eq} has no solution on $[0,\tau)$.
\end{Theorem}

Note that we are not asking $f$ to be below $\beta_0^{-\alpha}G$ in the whole interval $(0,\tau)$, but only at some specific time. Let us also remark that two of the main differences  with respect to the case $\alpha=1$ have already appeared: the fractional polynomials $P_k^\pm$  and the ordinary differential equation~\eqref{eq:edo}, that will be studied more profoundly for a particular $f$ later, in Section~\ref{sect-special_case}.

Next lemma, which is a corollary of Theorem~\ref{thm:no_existence_negative_kernel}, is needed for Theorem~\ref{thm:uniqueness_positive_kernel_or_0}, as \cite[Lemma 12]{Gilding-1993} is needed for the proofs of Lemma~13 and Theorem~16 in that paper.

\begin{Lemma}
Assume $\esssup\limits_{t\in(0,\tau)} g(t)\leq 0$ for some $\tau\in(0,\infty]$ and  $f(t)\leq kG(t)$ for all $t\in[0,\tau)$ and some $k<k_0$. Then $\tilde{x}(t;0)$ solves~\eqref{eq:formal_eq} on $[0,\tau)$ only if $\tilde{x}(t;0)=f(t)=G(t)=0$ for all $t\in[0,\tau)$.
\end{Lemma}
It is also needed to prove the final theorem of this subsection.
\begin{Theorem}
Assume that \eqref{eq:assumption.negative.kernel} and \eqref{eq:condition.positivity.G} hold, and $f(0)=0$.  If given any $t\in(0,\tau)$ there exist constants $\delta>0$ and $k<k_0$ such that
$$
\displaystyle f(s)-f(t)\leq k \{G^{\alpha +1}(s) - G^{\alpha +1}(t)\}^\frac{1}{\alpha + 1}\quad\text{for all }s\in(t, t+\delta),
$$
then either the equation~\eqref{eq:formal_eq} has no solution or it has a uncountable number of them.
\end{Theorem}

\subsection{Existence: non-negative kernels}\label{sect-Positive_Kernel}

This time we assume that the kernel is non-negative,
\begin{equation}
\label{eq:assumption.non-negative.kernel}
\essinf\limits_{t\in (0,\tau)} g(t)\geq 0\quad\text{for some }\tau\in (0,\infty].
\end{equation}
Our first existence result involves the function $G$ defined in~\eqref{eq:G}.

\begin{Theorem}\label{thm:positive_kernel_propiedad_acotacion}
Assume~\eqref{eq:assumption.non-negative.kernel}. If for some $-\infty < \gamma \leq \beta\leq \infty$,
$$
\gamma G(t)\leq f(t)\leq \beta G(t)\quad\text{for all }t\in[0,\tau),
$$
then equation~\eqref{eq:formal_eq} has a unique solution $\tilde{x}(t;0)$ on $[0,\tau)$ such that
$$
\beta_0^{-\alpha}G(t)\leq \tilde{x}(t;0) - f(t) \leq \gamma_0^{-\alpha}G(t),
$$
where $\gamma_0$ and $\beta_0$ are respectively the unique positive roots of the polynomials $P^-_\gamma$ and $P^-_\beta$.  If $\beta=\infty$ then we take $1/\beta_0^\alpha=0$.
\end{Theorem}
\begin{proof}

If $\beta = \infty$  we are only saying that $\tilde{x}(t;0)\geq f(t)$, which is trivial, since the kernel is non-negative. So, let us assume that $\beta<\infty$. To prove the right-hand estimate we argue by contradiction. Suppose that it is false, then there must exist an interval $[t_0,t_1]\subset [0,\tau)$ such that
\begin{equation}\label{eq:positive_kernel_edo}
Y(t_0) = \gamma_0^{-\alpha}G(t_0)
\quad\text{and } Y(t)> \gamma_0^{-\alpha}G(t)\text{ for all }t\in (t_0,t_1],
\end{equation}
where $Y(t):=\tilde{x}(t;0) - f(t)$.
Multiplying~\eqref{eq:formal_eq} by $Y'$
we get~\eqref{eq:edo}
for every $t\in(t_0,t_1)$, or, substituting $x=f+Y$ and  recalling that $f\geq \gamma G$,
$$
0\geq (\gamma G + Y)^\alpha Y' - G^\alpha G'.
$$
But, using~\eqref{eq:positive_kernel_edo}, we get
$$
0\geq (\gamma+\gamma_0^{-\alpha})^\alpha Y' -G'.
$$
Since $\gamma+\gamma_0^{-\alpha} = \gamma_0$, it is easy to get now, integrating from $t_0$ to $t_1$, a contradiction with~\eqref{eq:positive_kernel_edo} in $t_1$. The left-hand estimate follows analogously.
\end{proof}

The last result of this part can be proved as~\cite[Theorem 16]{Gilding-1993}.

\begin{Theorem}\label{thm:uniqueness_positive_kernel_or_0}
Suppose that $\essinf\limits_{0<t<\tau} g(t)\geq 0$ and that given any $t\in(0,\tau)$ there exist $\delta\in(0,t]$ and  $\beta<k_0$ such that
$$
\displaystyle f(s)-f(t)\leq \beta \{G^{\alpha +1}(t) - G^{\alpha +1}(s)\}^\frac{1}{\alpha + 1}\quad\text{for all }s\in(\delta, t).
$$
Then equation~\eqref{eq:formal_eq} has a unique solution on $[0,\tau)$. Moreover, if there exists a point $\tilde{s}\in (0,\tau)$ for which $\theta(\tilde{s})=0$, then $\theta(s)\equiv G(s)\equiv 0$ in $[0,\tilde{s}]$.
\end{Theorem}

\subsection{A special case}\label{sect-special_case}

We consider now equation~\eqref{eq:formal_eq} when
\begin{equation}
\label{eq:special_case}
\int_0^t g(s)\, {\rm d}s <0\text{ for all }t>0,\quad
f(t)=k G(t)\text{ for all }t\geq 0, \quad k\ge k_0.
\end{equation}

Again as in \cite{Gilding-Kersner-2004}, we are able to prove that our equation admits an uncountable number of solutions characterized by their behaviour as $t\downarrow 0$.

\begin{Theorem}\label{thm:special_case}
Suppose that~\eqref{eq:special_case} holds.
We denote the positive roots of $P^+_k$ by $z_1$ and $z_2$ when it has two, and by $z_0$ when it has only one.

\begin{itemize}
\item[\rm (i)]If $k>k_0$ then equation~\eqref{eq:formal_eq} admits the maximal solution
$$
\tilde{x}(t;0)=z_2G(t),
$$
and for each $\rho\in\mathbb{R}$ a unique solution $x_\rho$ such that
$$
x_\rho(t)=z_1G(t) + \rho G^{\gamma+1}(t) + \mathcal{O}(G^{2\gamma + 1}(t))\quad\text{as }t\downarrow 0,\quad \gamma = \frac{\alpha(k-z_1)}{z_1} - 1>0,
$$
with maximal interval of existence $[0,T_\rho)$, with $T_\rho=\infty$ if $\rho\geq 0$ and $T_\rho$ finite if $\rho<0$, and no other solutions. Moreover, if $G(t)\to\infty$ as $t\to \infty$ then for every $\rho>0$ there holds $x_\rho(t)\sim z_2G(t)$ as $t\to \infty$.

\item[\rm (ii)] If $k=k_0$ then equation~\eqref{eq:formal_eq} admits the maximal solution
$$
\tilde{x}(t;0)=z_0G(t),
$$
and for each $\rho\in\mathbb{R}$ a unique solution $x_\rho$ such that
$$
x_\rho(t)=(z_0G -\eta GL+ GL^2(\rho + \gamma J^{-1}))(t) + \mathcal{O}(GL^3J^{-2}(t))\quad\text{as }t\downarrow 0,
$$
where
$$
\eta=\frac{2z_0}{\alpha +1}=\frac{2k\alpha}{(\alpha +1)^2}\quad\text{and}\quad \gamma = \frac{2(\alpha+2)\eta}{3(\alpha + 1)},
$$
with maximal interval of existence $[0,T_\rho)$, with $T_\rho$ finite, and no other solutions.
\end{itemize}
In both cases the solutions $x_\rho$ are monotone with respect to the parameter $\rho$.
\end{Theorem}
\begin{proof}
Suppose that $x$ is a solution of our equation on $[0,\delta)$ for some $\delta$ positive. We define
$$
Y(t):=-\int_0^t\frac{g(s)}{x^\alpha (s)}\, {\rm d}s.
$$
We can multiply the equation by $Y'$ to obtain
$$
G'G^\alpha (kG-Y)^{1-\alpha} - Y'(kG-Y)=0.
$$
This equation is homogeneous, so defining $V=Y/G$ we arrive to
$$
\frac{G'}{G}=\frac{(k-V)^\alpha V'}{1-(k-V)^\alpha V}
$$
or, in other words,
\begin{equation}\label{eq:diferencial_V_G}
\frac{{\rm d}V}{{\rm d}G} = \frac{1-(k-V)^\alpha V}{G(k-V)^\alpha}.
\end{equation}
We have to study the trajectories that satisfy this equation.

Let us first study the case where $k>k_0$, which means that there are two different roots of $P_k^+$. Inspired by the results in~\cite{Gilding-1993}, we define, for a certain $\gamma\in\mathbb{R}$ to be defined later, a new independent variable $G^*= G^\gamma$, and~\eqref{eq:diferencial_V_G} becomes
$$
\frac{{\rm d}V}{{\rm d}G^*} = \frac{1-(k-V)^\alpha V}{\gamma G^*(k-V)^\alpha}.
$$
We have to study the trajectories $V(G^*)$ as $G^*\to 0$.

Clearly, we have that $V(0)=k-z_i$, $i=1$ or $2$. In fact, the only trajectory that can start from $k-z_{2}$ is the constant $V(G^*)=k-z_2$, which, recalling the definition of $V$, translates to a solution $x(t)=z_2G(t)$. More diverse is the family of trajectories sprouting from $k-z_1$. We would like to study the limit of ${\rm d}V/{\rm d}G^*$ as $G^*$ goes to 0.

Suppose that this limit exist, and define
$$
\rho = \lim\limits_{G^*\to 0} \frac{1-(k-V)^\alpha V}{\gamma G^*(k-V)^\alpha} = \frac{1}{\gamma z_1^\alpha} \lim\limits_{G^*\to 0} \frac{1-(k-V)^\alpha V}{G^*}.
$$
Then, applying L'H\^opital's rule, we can check that
$$
\rho = \frac{1}{\gamma z_1^\alpha} \lim\limits_{G^*\to 0} V' (\alpha V (k-V)^{\alpha-1} - (k-V)^\alpha)= \frac{\rho}{\gamma z_1^\alpha} (\alpha z_1^{\alpha -1}(k-z_1) - z_1^\alpha),
$$
meaning that $\rho$ is arbitrary whenever
$$
\gamma = \frac{\alpha(k-z_1)}{z_1} - 1.
$$
To see that this limit $\rho$ in fact exists we study the monotonicity of $(1-(k-V)^\alpha V)/G^*$ with respect to $G^*$. Indeed,
$$
\partial_{G^*} \left(\frac{1-(k-V)^\alpha V}{G^*}\right)= \frac{V'}{G^*} (\alpha V - (\gamma + 1)(k-V))= \frac{V'\alpha k}{z_1G^*} (V-(k-z_1)).
$$
Since the trajectories never cross the trajectory $V(G^*)=k-z_1$ (i.e. $V'(G^*)=0$), this means that the previous derivative has a constant sign, which gives monotonicity and existence of the limit $\rho$. The rest of the assertions when $k>k_0$ follow from the study of the phase-plane and the definition of $V$.

The case $k=k_0$, when we have only one root $z_0=\alpha k/(\alpha +1)$ of $P_k^+$, is more difficult, since $\gamma$ becomes 0. In this case we can say that the zero in the numerator $1-(k-V)^\alpha V$ is one order higher, so to compensate for it we would need a higher order zero in the denominator. We achieve this by considering ${\rm d}V/{\rm d}L$, recalling the definition of $L$. Summarizing,
\begin{equation}\label{eq:second_derivative_L}
\begin{array}{l}
\displaystyle
\frac{{\rm d}V}{{\rm d}L}= \frac{1-(k-V)^\alpha V}{L^2(k-V)^\alpha},\\[10pt]
\displaystyle\frac{{\rm d}^2V}{{\rm d}L^2}= \frac{V'}{k-V}\left( \frac{(\alpha +1)V - k}{L^2} - \frac{2(k-V)}{L} + \alpha V' \right).
\end{array}
\end{equation}

Now, from a study similar to the one before, we can see that
$$
\lim\limits_{L\to 0} V(L)=k-z_0\quad \text{and }\eta:=\lim\limits_{L\to 0} \frac{{\rm d}V}{{\rm d}L}= \frac{2z_0}{\alpha +1}
$$
for every trajectory. So, the distinction between them has to come from the second order derivative, but a problem appears since through L'H\^opital's rule we get
$$
\lim\limits_{L\to 0} \frac{{\rm d}^2V}{{\rm d}L^2}= \infty.
$$

Let us write then a possible expansion of the function $V(L)$ near $L=0$. We will omit the higher order terms  for the sake of simplicity, since they will not play an important role here. For a certain function $Q(L)$ we have, near zero, that
$$
V(L)\approx (k-z_0) + \eta L + Q(L)L^2,\quad V'(L)\approx \eta + 2LQ(L) + L^2Q'(L)
$$
and
$$
V''(L)\approx 2Q(L) + 4LQ'(L) + L^2Q''(L),
$$
where we impose that both $LQ(L)$ and $L^2Q'(L)$ go to 0 as $L$ goes to 0 in order for $V'(0)$ to be equal to $\eta$. If we substitute in~\eqref{eq:second_derivative_L} we get
$$
\begin{array}{lcc}
2Q(L) + 4LQ'(L) + L^2Q''(L) \approx \displaystyle\frac{V'}{k-V} \cdot\frac{(\alpha +1)((k-z_0) + \eta L + Q(L)L^2) - k}{L^2} \\[8pt]
\displaystyle\qquad - \frac{V'}{k-V} \left(\frac{2(z_0 - \eta L - Q(L)L^2)}{L} - \alpha V' \right).
\end{array}
$$
Since $(\alpha + 1)(k-z_0)-k=0$ and $(\alpha+1)\eta - 2z_0=0$ we can cancel the problematic terms and make $L$ approach 0 to get
$$
2Q(L) + 4LQ'(L) + L^2Q''(L) \approx \frac{\eta}{z_0} ((\alpha + 1)Q(L) + (\alpha + 2)\eta),
$$
which means that we can impose $4LQ'(L) + L^2Q''(L) = \eta^2(\alpha + 2)/z_0$. Thus $Q(L)\to \infty$ as $L\to 0$. The function that satisfies all these conditions is
$$
Q(L) = \rho - \frac{2(\alpha+2)\eta}{3(\alpha + 1)}\ln (L),
$$
where $\rho$ is once more an arbitrary real constant. Again an analysis of the  trajectories of the phase-plane and the definition of $V$ give the desired result.
\end{proof}

\noindent\emph{Remark. } If we let $\alpha \to 1$  we recover the result from Gilding in~\cite{Gilding-1993}.

%
%
%
%



\section{Travelling waves. General results}
\label{sect-tw.general.results} \setcounter{equation}{0}

It is time now to translate the results obtained in the previous section for the integral equation to the frame of travelling waves. Due to the possible degeneracies/singularities of the equation, in general we will need to deal with weak solutions. A function $V$ defined on an open real interval  $\Omega$ with values on $I$ is said to be a travelling-wave profile of equation~\eqref{eq:complete_equation} corresponding to the  speed $\sigma$ if $V\in C(\Omega),  |(a(V))'|^{p-2}(a(V))'$ and $c(V)\in L^1_{\rm loc}(\Omega)$ and
\begin{equation}\label{eq:weak_form_TW}
\int_\Omega |(a(V))'|^{p-2}(a(V))'\varphi' + b(V)\varphi' +\sigma V \varphi' - c(V)\varphi =0
\end{equation}
for every $\varphi \in C_0^\infty(\Omega)$. If $\Omega=\mathbb{R}$ then it is called a \textit{global} travelling-wave profile.

The constitutive functions $a, b, c: \bar I\to\mathbb{R}$ are only assumed to satisfy:
\begin{itemize}
\item[(i)] $a\in C(\bar I)\cap C^1(I)$, $a'(u)>0$ in $I$, and $a(0)=0$;
\item[(ii)] $b\in C(\bar I)\cap  C^1(I)$ and $b(0)=0$;
\item[(iii)] $c\in C(\bar I)$, with $c(0)=0$ and $c(l)=0$ if $l<\infty$.
\end{itemize}
It is important to note that much less can be asked from $c$, which needs not to be continuous for example, see~\cite{Gilding-Kersner-2004}, but we stick to the above assumptions for the sake of simplicity.

\subsection{Semiwavefronts}

If a profile is monotonic but is not defined everywhere but only in some interval $(\omega,\infty)$, and it still satisfies $V(\xi)\to \ell^+$  as $\xi\to \infty$ then it is called a \textit{semi-wavefront profile to} $\ell^+$. \emph{Semi-wavefronts profiles from} $\ell^-$ are defined in a similar way. The process of construction of wavefront profiles starts by considering first semi-wavefront profiles, and then extending them to the whole $\mathbb{R}$. A semi-wavefront profile that is not extendible to be global is called a \textit{strict semi-wavefront} profile. As an example we have
$$
V(\xi)=\left(\frac{\sigma(\alpha-m)}{m}\left|\frac{\alpha-m}{m}\right|^{p-2}\right)^{\frac{\alpha}{m-\alpha}}  \xi^\frac{1}{m-\alpha},\quad\xi>0,\qquad \alpha=1/(p-1),
$$
for any speed $\sigma> 0$, which is a semi-wavefront profile, not extendible to a wavefront, for the fast diffusion equation $u_t=\Delta_p u^m$ with $m-\alpha<0$.

We proceed to study the common results for travelling-wave profiles of equation~\eqref{eq:complete_equation} obtained from the study of the integral equation~\eqref{eq:complete_integral_equation}. Let us first see the direct applications of Theorem~\ref{thm:formal_comparison_maximal}.

\begin{Theorem}
If equation~\eqref{eq:complete_equation} has a semi-wavefront profile with speed $\sigma_0$ decreasing to 0 then it has such a solution for every wave speed $\sigma\geq\sigma_0$. Consider also equation~\eqref{eq:complete_equation} with two sets of coefficients $a_i,b_i$ and $c_i$ for $i=1,2$ on some interval $I$ and two real parameters $\sigma_1$ and $\sigma_2$.
\begin{itemize}
\item[\rm (a)] Suppose that the function $u\mapsto \sigma_2 u + b_2(u)-\sigma_1 u - b_1(u)$ is nondecreasing on $ I$ and $(c_2a_2')(u)\leq (c_1a_1')(u)$ for all $u\in(0,\ell)$.
\item[\rm (b)] Else, suppose that $\sigma_2 u + b_2(u)\geq\sigma_1 u + b_1(u)$ and $\max\{0,(c_2a_2')(u)\}\leq (c_1a_1')(u)$ for all $u\in(0,\ell)$.
\end{itemize}
Then in both cases, if equation~\eqref{eq:complete_equation} with $i=1$ admits a semi-wavefront profile decreasing to 0 with speed $\sigma_1$, so does it with $i=2$ and speed $\sigma_2$.
\end{Theorem}

This result implies the existence of a wave speed $\sigma^*\in[-\infty,\infty]$ such that \eqref{eq:complete_equation} admits semi-wavefront profiles for  $\sigma >\sigma^*$  but not for $\sigma<\sigma^*$. Due to this, the speed $\sigma^*$ is called the \textit{critical speed} of the equation.

\noindent\emph{Remark. } The semi-wavefront profiles corresponding to speeds $\sigma>\sigma_0$ do not necessarily have the same properties as the one with speed $\sigma_0$. For instance, the  semi-wavefront profile for the speed $\sigma_0$ may vanish for large values of $\xi$ while the  ones for speeds $\sigma>\sigma_0$ are positive everywhere.

\medskip

Let us now have some words about the number of semi-wavefront profiles for a fixed wave speed $\sigma$. We say that the equation~\eqref{eq:complete_equation} admits a one-parameter family of distinct semi-wavefront profiles  decreasing to 0 with speed $\sigma$ when there exists a continuous order-preserving bijective mapping from the interval~$(0,1]$ onto the set of such solutions.

\begin{Theorem}\label{thm:number_solution_sink_general_eq}
{\rm (i)}
If $c<0 $ in $(0,\ell)$, then equation~\eqref{eq:complete_equation} has at most one distinct semi-wavefront profile decreasing to 0 with wave speed $\sigma$.

\noindent {\rm (ii)} If $c>0$ in $(0,\ell)$,  then equation~\eqref{eq:complete_equation} has either a one-parameter family of distinct semi-wavefront profiles decreasing to~0 with wave speed $\sigma$  or no such solution for that speed.
\end{Theorem}

The first assertion is a consequence of Theorem~\ref{thm:formal_uniqueness} and  Corollary~\ref{thm:es_el_2.25_del_libro}. The second one is more complicated and requires Corollary~\ref{thm:es_el_2.25_del_libro} and some extra work, namely two lemmata that are similar to \cite[Lemmas 4.7 and 4.8]{Gilding-Kersner-2004},  whose proof we omit for the sake of brevity.

\subsection{Wavefronts}

Let us start this subsection by presenting three equivalence results.
\begin{Theorem}
Suppose $\ell<\infty$. If  one of the equations
$$
u_t=\Delta_p (a(u))+  \nabla(b(u)) + c(u),\qquad
u_t=\Delta_p u + \nabla(b(u)) + c(u)a'(u)
$$
admits a wavefront profile from $\ell$ to 0 with speed $\sigma$, then they both do. Moreover if $a\in C^1(I)$, there is an explicit transformation from a  wavefront profile $V_1$ for the first one onto a wavefront profile $V_2$ for the second one given by
$$
V_2(\xi)=V_1(\Psi(\xi)), \quad\xi\in\mathbb{R},\qquad \text{where }
\Psi(\xi)=\int_0^\xi a'(V_2(\eta))\, {\rm d}\eta.
$$
\end{Theorem}

This equivalence  was first studied by Engler in~\cite{Engler-1985} when $\alpha=1$ and $a(u)=u^m$.  Gilding and Kersner extended the result in~\cite{Gilding-Kersner-2004} to more general nonlinearities $a$. Their proof works with only trivial changes for $\alpha\neq1$.

\begin{Theorem}\label{thm:correspondance_u_l-u}
Suppose $\ell<\infty$. Equation~\eqref{eq:complete_equation} has a wavefront profile from $\ell$ to~0 with speed $\sigma$ if and only if the equation
\begin{equation}\label{eq:ecuacion_contraria}
\begin{array}{c}
u_t=\Delta_p(\tilde{a}(u)) + \nabla(\tilde{b}(u)) + \tilde{c}(u), \text{ where } \\[8pt]
\tilde{a}(u)= a(\ell)-a(\ell-u),\quad \tilde{b}(u) = b(\ell-u)-b(\ell),\quad\text{and }\tilde{c}(u) = -c(\ell-u),
\end{array}
\end{equation}
has a wavefront profile from  0 to $\ell$ with speed $-\sigma$. In both cases the number of distinct wavefront profiles is the same.
\end{Theorem}

One can check that if $u$ is a wavefront solution of our original equation that connects $\ell$ with~0 then $v(x,t)=\ell-u(-x,t)$ is a solution of~\eqref{eq:ecuacion_contraria} that also connects $\ell$ with~0. In terms of the integral equation, this result is a consequence of the following theorem.

\begin{Theorem}\label{thm:correspondance_theta_Theta}
Suppose that $\ell<\infty$. Then the following statements are equivalent: equation~\eqref{eq:complete_integral_equation} has a solution $\theta$ on $\bar{I}$ with $\theta(l)=0$; equation
\begin{equation}\label{eq:ecuacion_integral_contraria}
\Theta(s)=-\sigma s + \tilde{b}(s) - \int_0^s\frac{\tilde{c}(r)\tilde{a}'(r)}{\Theta^\alpha(r)} \, {\rm d}r
\end{equation}
has a solution $\Theta$ on $\bar{I}$ with $\Theta(\ell)=0$; and, equations~\eqref{eq:complete_integral_equation} and~\eqref{eq:ecuacion_integral_contraria} both have solutions on $\bar{I}$. In the same way, these statements are also equivalent: equation~\eqref{eq:complete_integral_equation} has a solution $\theta$ satisfying the integrability condition on $\bar{I}$ with $\theta(l)=0$; equation~\eqref{eq:ecuacion_integral_contraria} has a solution $\Theta$ satisfying the integrability condition on $\bar{I}$ with $\Theta(l)=0$; and, equations~\eqref{eq:complete_integral_equation} and~\eqref{eq:ecuacion_integral_contraria} both have solutions satisfying the integrability condition on $\bar{I}$.
\end{Theorem}

These results give us the next one, about the admissible wave speeds; see~\cite[Theorem~8.1]{Gilding-Kersner-2004}.

\begin{Theorem}
Suppose that $\ell<\infty$. If equation~\eqref{eq:complete_equation} has wavefront solutions from $\ell$ to 0 with speed $\sigma_1$ and with speed $\sigma_2>\sigma_1$, the same is true for all $\sigma\in[\sigma_1, \sigma_2]$.
\end{Theorem}

The result can be improved if we impose additional assumptions on the reaction term $c$.
\begin{Theorem}
Suppose that $\ell<\infty$.
\begin{itemize}
\item[\rm (i)] If $c<0$ in $(0,\ell)$, the set $S$ of wave speeds for which~\eqref{eq:complete_equation} has a wavefront solution from $\ell$ to 0 is either empty or there exists $\sigma_0$ such that $S=(-\infty, \sigma_0]$.
\item[\rm (ii)] If $c>0$ in $(0,\ell)$, either $S=\emptyset$ or $S=[\sigma_0,\infty)$ for some value $\sigma_0$.
\item[\rm (iii)] If $c\leq 0$ in $(0,\ell)$, the set $S$ is either empty, contains a single value or is an interval which is bounded above and contains its right endpoint.
\item[\rm (iv)] If $c\geq 0$ in $(0,\ell)$, the set $S$ is either empty, contains a single value or is an interval which is bounded below and contains its left endpoint.
\end{itemize}
\end{Theorem}
This corresponds to~\cite[Theorem 8.3]{Gilding-Kersner-2004}. The proof, though lengthy, is again similar to the one in it.

Let us now have some words about the number of solutions for a fixed speed. This result is a consequence of Theorems~\ref{thm:number_solution_sink_general_eq} and~\ref{thm:correspondance_u_l-u} and corresponds to~\cite[Theorem 8.7]{Gilding-Kersner-2004}.

\begin{Theorem}
Suppose that $\ell<\infty$ and let $\sigma$ be a fixed wave speed. Then equation~\eqref{eq:complete_equation} has at most one distinct wavefront solution from $\ell$ to 0 with speed $\sigma$ whenever $c<0$ in $(0,\ell)$ or $c>0$  in $(0,\ell)$.
\end{Theorem}

Finally we focus on the support of the wavefront profiles. It is well known that some combinations of filtration and reaction nonlinearities $a$ and $c$  provoke the appearance of free boundaries in the solutions. Thus, for wavefront profiles from $\ell$ to 0 one or both of the following properties may hold:
\begin{eqnarray}
\label{eq:suport_bouded_above}
V(\xi)&\equiv 0\quad\text{for all } \xi\geq \xi^*\quad\text{for some }\xi^*\in\mathbb{R},\\
\label{eq:suport_bouded_below}
V(\xi)&\equiv \ell\quad\text{for all } \xi\leq \xi_*\quad\text{for some }\xi_*\in\mathbb{R}.
\end{eqnarray}
The same is true for semi-wavefront solutions. For the sake of brevity, we only present the result for wavefronts, noticing that it can be adapted to semi-wavefronts; see~\cite[Theorems 2.30, 2.34 and 2.38]{Gilding-Kersner-2004} for more information about this.

\begin{Theorem}\label{thm:bounded_support_in_general} Equation~\eqref{eq:complete_equation} has a wavefront profile $V$ with speed $\sigma$ such that:
\begin{itemize}
\item[\rm (i)] Properties~\eqref{eq:suport_bouded_above} and~\eqref{eq:suport_bouded_below}
hold simultaneously if and only if equation~\eqref{eq:complete_integral_equation} has a solution $\theta$ on $\bar{I}$ such that $\theta(\ell)=0$ and
$$
\int_0^\ell \frac{a'(r)}{\theta^\alpha(r)}\, {\rm d}r <  \infty.
$$
\item[\rm (ii)]
Property~\eqref{eq:suport_bouded_above} and
\begin{equation}\label{eq:suport_not_bouded_below}
V(\xi)< \ell\quad\text{for all } \xi\in\mathbb{R}
\end{equation}
hold simultaneously
if and only if equation~\eqref{eq:complete_integral_equation} has a solution $\theta$ on $\bar{I}$ such that
$$
\int_0^s \frac{a'(r)}{\theta^\alpha(r)}\, {\rm d}r <  \infty\quad \text{for all }s\in(0,\ell),
$$
and, either
\begin{equation}\label{eq:integral_entera_de_0_a_l_para_soporte}
\int_0^\ell \frac{a'(r)}{\theta^\alpha(r)}\, {\rm d}r =  \infty
\end{equation}
or $\theta(s_i)=0$ for a sequence of values $\{s_i\}_{i=1}^\infty\subset I$ such that $s_i\to \ell$ as $i\to \infty$.
\item[\rm (iii)] Property~\eqref{eq:suport_bouded_below}  and
\begin{equation}\label{eq:suport_not_bouded_above}
V(\xi)>0\quad\text{for all } \xi\in\mathbb{R}
\end{equation}
hold simultaneously if and only if~\eqref{eq:complete_integral_equation} has a solution $\theta$ on $\bar{I}$ such that $\theta(l)=0$,
$$
\int_s^\ell \frac{a'(r)}{\theta^\alpha(r)}\, {\rm d}r <  \infty\quad \text{for all }s\in(0,\ell),
$$
and, either~\eqref{eq:integral_entera_de_0_a_l_para_soporte} is satisfied or $\theta(s_i)=0$ for a sequence of values $\{s_i\}_{i=1}^\infty\subset I$ such that $s_i\to 0$ as $i\to \infty$.
\item[\rm (iv)] Properties~\eqref{eq:suport_not_bouded_below} and~\eqref{eq:suport_not_bouded_above} hold simultaneously if and only if equation~\eqref{eq:complete_integral_equation} has a solution $\theta$ satisfying the integrability condition on $\bar{I}$ such that
$$
\int_0^s \frac{a'(r)}{\theta^\alpha(r)}\, {\rm d}r = \infty\quad \text{for all }s\in(0,\ell),
$$
or $\theta(s_i)=0$ for a sequence of values $\{s_i\}_{i=1}^\infty\subset I$ such that $s_i\to 0$ as $i\to \infty$, and
$$
\int_s^\ell \frac{a'(r)}{\theta^\alpha(r)}\, {\rm d}r = \infty\quad \text{for all }s\in(0,\ell),
$$
or $\theta(s_i)=0$ for a sequence of values $\{s_i\}_{i=1}^\infty\subset I$ such that $s_i\to \ell$ as $i\to \infty$.
\end{itemize}
\end{Theorem}

\section{Reaction-diffusion}
\label{sect-travelling_waves_reaction_diffusion} \setcounter{equation}{0}

In this section we will study in more detail travelling waves for the equation~\eqref{eq:complete_equation} when there is no convection, $b\equiv 0$, which leaves us with the reaction-diffusion equation
\begin{equation}\label{eq:differential_reaction_diffusion}
u_t=\Delta_p (a(u)) + c(u).
\end{equation}
For an equation of this class the integral equation becomes
\begin{equation}\label{eq:integral_reaction_diffusion}
\theta(s) = \sigma s  - \int_0^s \frac{a^\prime (r) c(r)}{\theta^\alpha (r)}\, {\rm d}r,\qquad \alpha>0.
\end{equation}

\noindent\emph{Remark. }  If $c\equiv 0$ ($b\not\equiv0$), the integral equation~\eqref{eq:complete_integral_equation} becomes $\theta(s)=\sigma s+ b(s)$, and the search for nonnegative solutions in an interval $(0,\delta)$ is much easier. Nevertheless, one has to treat the possibility of having $\theta(s^*)=0$ for some $s^*>0$ more carefully than in the case $p=2$, since on the one hand this could happen without breaking the integrability condition~\eqref{eq:integrability_condition} and on the other hand there is no reaction to ``compensate" the degeneracy. Anyway, the analysis of this case is not too hard, and is analogous to the one in~\cite{Gilding-Kersner-2004}, hence we skip it.

\subsection{Semiwavefronts}

We consider here the  two cases in which the reaction term $c$ has a definite sign, either negative (sink) or positive (source). These two cases will later be combined to produce wavefront profiles for reaction nonlinearities with one sign change.

\subsubsection{Sink term}\label{sect-sink_term_travelling_waves}

We start by considering the case in which $c<0$ in $(0,\ell)$. We start with a uniqueness and existence result, which is a consequence of Theorems~\ref{thm:uniqueness_positive_kernel_or_0} and~\ref{thm:number_solution_sink_general_eq}.

\begin{Theorem}\label{thm:existance_negative_reaction}
Suppose that $c<0$ in $(0,\ell)$. Then for every wave speed $\sigma$ equation~\eqref{eq:differential_reaction_diffusion} has exactly one distinct semi-wavefront profile decreasing to 0.
\end{Theorem}

We worry now about the support of these solutions. We start with a lemma which is analogous to \cite[Lemma 6.3]{Gilding-Kersner-2004}, that gives us some barriers for the flux.

\begin{Lemma}
Suppose that $c\leq 0$ in $(0,\delta)$ for some $\delta\in(0,\ell)$, and that equation~\eqref{eq:integral_reaction_diffusion} has a unique solution $\theta$ on $[0,\delta)$.
\begin{itemize}
\item[\rm (i)] If $\sigma=0$ then $\theta(s)=G(s)$ for all $s\in[0,\delta)$.

\item[\rm (ii)] If $\sigma>0$ then
$$
\min\{\sigma, 1\}\max\{s,G(s)\}\leq \theta(s)\leq (\sigma+1)\max\{s,G(s)\}\quad\text{for all }s\in(0,\delta).
$$
Moreover, if $G(s)/s\to \mu$ as $s\to 0$ for some $0\leq \mu\leq \infty$ then
$\theta(s)/s\to \mu_0$ as $s\downarrow 0$,
where $\mu_0$ is positive constant depending on $\mu$. Finally, if $\mu>0$ then
$$
\frac{\theta(s)}{G(s)}\to \frac{\mu_0}{\mu}\quad\text{as }s\downarrow 0,\quad \frac{\mu_0}{\mu}\to 1 \quad\text{as }\mu\to\infty.
$$

\item[\rm (iii)] If $\sigma <0$ and $|c(s)a'(s)|\geq AG^\alpha(s)$ for almost all $s\in(0,\delta)$ for some $A>0$. then
$$
\theta(s)\geq z_0  G(s)\quad\text{for all }s\in(0,\delta),
$$
where $z_0$ is the unique positive root of  $P^-_{\sigma/A}$.

\item[\rm (iv)] If $\sigma<0$, $ca'$ is absolutely continuous on $[0,\delta)$,  $(|ca'|^{1/\alpha})'(s)\leq A^{1/\alpha}$  for almost all $s\in(0,\delta)$ for some $A\geq 0$, and $c(0)a'(0)=0$, then
$$
\theta(s)\geq z_1 |c(s)a'(s)|^{1/\alpha}\quad\text{for all }s\in(0,\delta),
$$
where $z_1$ is the unique positive root of  $P^*_A(z):=A^{1/\alpha}z^{\alpha+1} + |\sigma|z^\alpha -1$.

\item[\rm (v)] If $\sigma<0$, $ca'$ is absolutely continuous on $[0,\delta)$,  $(|ca'|^{1/\alpha})'(s)\geq B^{1/\alpha}$ for almost all $s\in(0,\delta)$ for some $B\geq 0$, and $c(0)a'(0)=0$, then
$$
\theta(s)\leq z_2 |c(s)a'(s)|^{1/\alpha}\quad\text{for all }s\in(0,\delta),
$$
where $z_2$ is the unique positive root of the fractional polynomial $P^*_B$.
\end{itemize}
\end{Lemma}
\begin{proof} (i) Taking $\gamma=\beta=0$ in Theorem~\ref{thm:positive_kernel_propiedad_acotacion} yields this result.

\noindent (ii) Taking $\gamma=0$ and $\beta=\infty$ in Theorem~\ref{thm:positive_kernel_propiedad_acotacion} gives $\sigma s\leq \theta(s)\leq \sigma s + G(s)$ for all $0\leq s < \delta$. However, Theorems~\ref{thm:formal_uniqueness} and~\ref{thm:formal_comparison_maximal} and part (i) of this lemma also provide $\theta(s)\geq G(s)$ for all $s\in [0,\delta)$. These two inequalities combined provide the first assertion.

To prove the next part we take $\delta^*\in(0,\delta)$ and define $A:=\inf\limits_{s\in(0,\delta^*)}(G(s)/s)$ and $B:=\sup\limits_{s\in(0,\delta^*)}(G(s)/s)$. Taking now $\gamma=\sigma/B$ and $\beta=\sigma/A$ in Theorem~\ref{thm:positive_kernel_propiedad_acotacion} gives
$$
\sigma s  + \beta_0 G(s)\leq \theta(s)\leq \sigma s  + \gamma_0
G(s)$$
for a certain couple of positive values $\beta_0,\gamma_0$. Passing to the limit as $\delta^*\downarrow 0$ gives the other two assertions, since $\beta_0$ and $\gamma_0$ both must converge to a certain $\mu_0$ dependant on $\mu$, whose behaviour as $\mu\to\infty$ comes from Theorem~\ref{thm:positive_kernel_propiedad_acotacion}.

\noindent (iii) One can check that $\theta_1(s):=z_0G(s)$ is a solution of our equation on $[0,\delta)$ with $\sigma_1 s + b_1(s):= (\sigma/A)G(s)$, and also that $s\mapsto \sigma s - \sigma_1 s - b_1(s)$ is nondecreasing on $[0\,\delta)$. The assertion comes then from Theorems~\ref{thm:formal_uniqueness} and~\ref{thm:formal_comparison_maximal}, because they mean that the unique solution satisfies $\theta\geq \theta_1$ on $[0,\delta)$.

\noindent (iv) This part is similar to the previous one considering $\theta_1(s):=z_1|c(s)a'(s)|^{1/\alpha}$ and $\sigma_1 s+b_1(s):=\theta_1(s) - 1/kz_1^\alpha s$.

\noindent (v) The function $\theta_2(s):=z_2|c(s)a'(s)|^{1/\alpha}$ is a solution of our equation on $[0,\delta)$ with $\sigma_2 s + b_2(s):= \theta_2(s) - 1/z_2^\alpha s$, and $s\mapsto \sigma_2 s + b_2(s) - \sigma s$ is nondecreasing in $[0,\delta)$. We finish as in (iii).
\end{proof}

As a consequence of this lemma and Theorem~\ref{thm:bounded_support_in_general} we obtain conditions guaranteing either the positivity of the wavefront profile or, on the contrary, that its support is bounded from above.
\begin{Theorem}
Suppose that the conditions of Theorem~\ref{thm:existance_negative_reaction} hold. Fix $\delta\in(0,\ell)$.
\begin{itemize}
\item[\rm (i)] If
$$
\int_0^\delta \frac{a'(s)}{(\max\{s,G(s)\})^\alpha}\, {\rm d}s =\infty,
$$
then every semi-wavefront profile decreasing to 0 is positive everywhere in its domain of definition.

\item[\rm (ii)]  If
$$
\int_0^\delta \frac{a'(s)}{(\max\{s,G(s)\})^\alpha}\, {\rm d}s <\infty\text{ and }
\int_0^\delta \frac{a'(s)}{G(s)^\alpha}\, {\rm d}s =\infty,
$$
then every semi-wavefront profile decreasing to 0 with speed $\sigma\leq 0$ is positive everywhere in its domain of definition and the ones with $\sigma>0$ have a support bounded from above.

\item[\rm (iii)] If
$$
\int_0^\delta \frac{a'(s)}{G(s)^\alpha}\, {\rm d}s <\infty,
$$
then every semi-wavefront profile decreasing to 0 with speed $\sigma\geq 0$  has a support bounded from above. Moreover, if $ca'$ is absolutely continuous in $[0,\delta)$, $(ca')(0)=0$,
$$
\lim \esssup\limits_{s\to 0} ((ca')^{1/\alpha})'(s)\leq 0\text{ and }
\int_0^\delta \frac{1}{|c(s)|}\, {\rm d}s =\infty,
$$
then every solution of this type with speed $\sigma<0$ is positive everywhere in its domain of definition. On the contrary, if $\lim\sup\limits_{s\to 0}G(s)/s=\lim \inf\limits_{s\to 0} G(s)/s>0$, or if $\lim \essinf\limits_{s\to 0} ((ca')^{1/\alpha})'(s)/G(s)> 0$, or if $ca'$ is absolutely continuous on $[0,\delta)$, $(ca')(0)=0$, $\lim \essinf\limits_{s\to 0} ((ca')^{1/\alpha})'(s)>-\infty$ and
$$
\int_0^\delta \frac{1}{|c(s)|}\, {\rm d}s <\infty
$$
then every solution of this type with speed $\sigma<0$ has its support bounded above.
\end{itemize}
\end{Theorem}

\subsubsection{Source term}\label{sect-source_term_travelling_waves}

The case $c>0$ in $(0,\ell)$ is more complicated. For an instance, our equation may not have a semi-wavefront profile. To study it, we define the quantities
$$
\lambda_1:=\limsup\limits_{s\downarrow 0} \left\{ \frac{1}{s}\int_0^s\frac{c(r)a'(r)}{r^\alpha} \, {\rm d}r \right\},\quad \Lambda_1(\delta):=\sup\limits_{0<s<\delta}\left\{ \frac{1}{s}\int_0^s\frac{c(r)a'(r)}{r^\alpha} \, {\rm d}r \right\},
$$
$$
\lambda_0:=\liminf\limits_{s\downarrow 0} \left\{ \frac{1}{s}\int_0^s\frac{c(r)a'(r)}{r^\alpha} \, {\rm d}r \right\}\text{ and } \Lambda_0(\delta):=\inf\limits_{0<s<\delta}\left\{ \frac{1}{s}\int_0^s\frac{c(r)a'(r)}{r^\alpha} \, {\rm d}r \right\}.
$$
We also define
$\sigma_{\lambda}:= k_0\lambda^{\frac{1}{\alpha+1}}$.

The following lemma, similar to \cite[Lemma 6.6]{Gilding-Kersner-2004}, will be needed for our existence result. We will write simply $\Lambda_{0}$ or $\Lambda_{1}$ when the dependence on $\delta$ is clear.
\begin{Lemma}\label{lemma:existence_reaction_difusion_c_positive}
Suppose that $c>0$ in $(0,\delta)$ for some $\delta\in (0,\ell)$.
\begin{itemize}
\item[\rm (i)] If $\sigma\leq 0$ or $\sigma <\sigma_{\Lambda_0}$, the integral equation~\eqref{eq:integral_reaction_diffusion} has no solution on $[0,\delta)$. Furthermore, any solution $\theta$ of this equation satisfies $\theta(s)\leq \mu s$ for all $s\in(0,\delta)$, where $\mu$ is the biggest root of $Q_{\Lambda_0}$, where
    $$
    Q_\gamma(z):=z^{\alpha +1} - \sigma z^\alpha + \gamma.
    $$

\item[\rm (ii)] If $\Lambda_1=\infty$ or if
$$
(\alpha+1)\Lambda_1 < \mu^\alpha [(\alpha+1)\sigma + \mu((1-\alpha)\alpha^{\frac{\alpha +1}{\alpha-1}} -1)],
$$
then~\eqref{eq:integral_reaction_diffusion} has no solution on $[0,\delta)$. We recall that $\mu$ depends only on $\sigma$ and $\Lambda_0$.

\item[\rm (iii)]  If $\lambda_1 < \Lambda_1<\infty$ and
$$
(\alpha+1)\Lambda_1 \leq \mu^\alpha [(\alpha+1)\sigma + \mu((1-\alpha)\alpha^{\frac{\alpha +1}{\alpha-1}} -1)]
$$
equation~\eqref{eq:integral_reaction_diffusion} has no solution on $[0,\delta)$.

\item[\rm (iv)] If $\sigma \geq \sigma_{\Lambda_1}$, equation~\eqref{eq:integral_reaction_diffusion} has a solution $\theta$ on $[0,\delta)$ such that $\theta (s)\leq \nu s$ for all $s\in[0,\delta)$, where $\nu$ is the biggest root of $Q_{\Lambda_1}$.
\end{itemize}
\end{Lemma}

\begin{proof}
\noindent (i) This part is similar to the one in \cite[Lemma 6.6]{Gilding-Kersner-2004}, after performing an analysis of the polynomial $Q_{\Lambda_0}$ analogous to the one for the polynomials $P_\gamma^\pm$ done in Section~\ref{sect-integral_equation}. Let us remark that this analysis yields $\mu\leq \sigma$, which will be important later.

\noindent(ii) We differentiate in~\eqref{eq:integral_reaction_diffusion} and multiply by $(\alpha +1)\theta^\alpha/s^\alpha$ to obtain
$$
\begin{aligned}
\frac{\theta^{\alpha+1}}{s^\alpha} &+ \frac{((\alpha+1)\sigma s - \mu s - \alpha \theta)(\mu^{\alpha}s^{\alpha} - \theta^\alpha)}{s^{\alpha +1}} + \frac{(\alpha+1)ca'}{s^\alpha}
\\
&=\mu^\alpha(\alpha +1)\sigma - \mu^{\alpha +1} + \frac{\alpha\theta^\alpha\mu}{s^\alpha} - \frac{\theta\mu^\alpha}{s}.
\end{aligned}
$$
On the other hand, the function $f(w)=\alpha\mu w^\alpha - \mu^\alpha w$ attains its maximum for $w\geq 0$ at $w_{max}= \mu \alpha^{2/(1-\alpha)}$, and therefore
$$
f(\theta/s)= \frac{\alpha\theta^\alpha\mu}{s^\alpha} - \frac{\theta\mu^\alpha}{s} \leq \mu^{\alpha +1} \alpha^{2/(1-\alpha)} \left(\frac{1}{\alpha} -1\right).
$$
Knowing this, recalling that $\theta(s) \leq \mu s \leq \sigma s$, integrating from 0 to $s\in (0,\delta)$ and dividing by $s$ we obtain our result.

\noindent (iii) If $\lambda_1 < \Lambda_1<\infty$ then there must exist a value $s\in(0,\delta)$ for which
$$
\Lambda_1=\frac{1}{s}\int_0^s\frac{c(s)a'(s)}{s^\alpha}\, {\rm d}s
$$
and this allows us to obtain a strict inequality in the previous part.

\noindent (iv) This part is similar to the one in \cite[Lemma 6.6]{Gilding-Kersner-2004}, and again uses our knowledge of the function $Q_{\Lambda_1}(z)$, the comparison results and the fact that $\theta_1(s):=\nu s$ is a solution of a similar equation with different coefficients.
\end{proof}

The main result about existence of solutions is the following.
\begin{Corollary}\label{thm:existence_positive_reaction}
Suppose that $c>0$ in $(0,\ell)$. Then for every wave speed $\sigma$ equation~\eqref{eq:differential_reaction_diffusion} has a one-parameter family of distinct semi-wavefront profiles decreasing to 0 or no such solution. Moreover:
\begin{itemize}
\item[\rm (i)] When $\lambda_1=\infty$ the equation has no such solution for all $\sigma$.
\item[\rm (ii)] When $\lambda_1\in(0,\infty)$, there exists a value $\sigma^*>0$ such that the equation has a one-parameter family of such solutions for all $\sigma>\sigma^*$ and no such solutions for $\sigma<\sigma^*$. The critical wave speed satisfies
$\sigma^*\in [\sigma_{\lambda_0},\sigma_{\lambda_1}]$.
Furthermore, if $\lambda_0=\Lambda_1(\delta)$ for some $\delta\in(0,\ell)$, the equation has a one parameter family of distinct semi-wavefront profiles decreasing to 0 with the critical speed~$\sigma^*$.
\item[\rm (iii)] When $\lambda_1=0$ the equation has a one-parameter family of such solutions for all $\sigma>0$ and no such solution for all $\sigma\leq 0$.
\end{itemize}
\end{Corollary}

This theorem is a consequence of Lemma~\ref{lemma:existence_reaction_difusion_c_positive}, once one notes that $\delta$ can be taken arbitrarily small.

Let us now present two lemmas useful for the study of the support of the wavefront profiles similar to \cite[Lemmata 6.7 and 6.8]{Gilding-Kersner-2004}. For the proof of the first one, which does not require $b\equiv0$, we refer to this book.

\begin{Lemma}\label{lemma:comparison_TW_reaction_diffusion}
Consider equation~\eqref{eq:complete_integral_equation} with two different wave speeds $\sigma_{i}$, sets of coefficients $a_i, b_i$ and $c_i$, and solutions $\theta_i$ on $[0,\delta)$,  $\delta\in(0,\ell]$. Suppose that the function $s \mapsto (\sigma_2 -\sigma_1)s + (b_2(s)-b_1(s))$ is nondecreasing on $[0,\delta)$, $c_2a'_2 \leq c_1a'_1$ almost everywhere in $(0,\delta)$ and $c_1(s)\neq 0$ or $c_2(s)\neq 0$ for all $0<s<\delta$. Then, either $\theta_2\geq \theta_1$ on $[0,\delta)$ or there exists  $\delta^*\in (0,\delta)$ such that $\theta_2 < \theta_1$ in $(0,\delta^*)$.
\end{Lemma}

\begin{Lemma}\label{lemma:behaviour_near_0_reaction_diffusion_c_positive}
Let the assumptions of Lemma~\ref{lemma:existence_reaction_difusion_c_positive} hold.
\begin{itemize}
\item[\rm (i)]  Let $\gamma^*$ be the smallest positive root of $Q_{\Lambda_1}$. If $\sigma > \sigma_{\Lambda_1}$, then given any $\gamma >\gamma^*$ equation~\eqref{eq:integral_reaction_diffusion} has at most one solution $\theta$ on $[0,\delta)$ such that
$$
\theta(s)\geq \gamma s\quad\text{for all }s\in(0,\delta).
$$

\item[\rm (ii)] If  $c(s)a'(s) \geq As^\alpha$ for almost all $s\in(0,\delta)$ for some $A>0$,  then given any solution $\theta$ of~\eqref{eq:integral_reaction_diffusion} on $[0,\delta)$ necessarily
$$
\sigma\geq\sigma_A\text{ and }
\liminf\limits_{s\downarrow 0}\frac{\theta(s)}{s}\geq z_1
$$
where $z_1$ is the smallest positive root of $Q_A$.

\item[\rm (iii)] If $\sigma > \sigma_{\Lambda_1}$ and $c(s)a'(s) \leq Bs^\alpha$ for almost all $s\in (0,\delta)$  for some $B\in(\Lambda_1,(\sigma/k_0)^{\alpha +1}]$,  then given any solution $\theta$ of~\eqref{eq:integral_reaction_diffusion} on $[0,\delta)$ other than the maximal solution necessarily
$$
\limsup\limits_{s\downarrow 0} \frac{\theta(s)}{s}\leq \beta_1
$$
where $\beta_1$ is the smallest positive root of  $Q_B$.

\item[\rm (iv)] If $ca'$ is absolutely continuous on $[0,\delta)$, $(ca')(0)=0$ and $((ca')^{1/\alpha})'(s)\geq A^{1/\alpha}$ for almost all $s\in(0,\delta)$ for some $A$, then given any solution $\theta$ of~\eqref{eq:integral_reaction_diffusion} on $[0,\delta)$ necessarily
$$
\sigma\geq \sigma_A\text{ and }\liminf\limits_{s\downarrow 0}\frac{\theta(s)}{(c(s)a'(s))^{1/\alpha}}\geq \frac{z_1}{A^{1/\alpha}}.
$$

\item[\rm (v)] If $\sigma >\sigma_{\Lambda_1}$, $ca'$ is absolutely continuous on $[0,\delta)$,  $((ca')^{1/\alpha})'(s)\leq B^{1/\alpha}$ for almost all $s\in(0,\delta)$ for some $B\in (\Lambda_1,(\sigma/k_0)^{\alpha +1}]$, and $(ca')(0)=0$,  then given any solution $\theta$ of~\eqref{eq:integral_reaction_diffusion} on $[0,\delta)$ other than the maximal solution, necessarily
$$
\limsup\limits_{s\downarrow 0}\frac{\theta(s)}{(c(s)a'(s))^{1/\alpha}}\leq \frac{\beta_1}{B^{1/\alpha}}.
$$
\end{itemize}
\end{Lemma}
\begin{proof} (i)
Let $X$ denote the set of real functions $\psi$ defined on $[0,\delta]$ such that $\nu~\leq~\psi\leq \sigma$, where  $\nu$ is again the biggest root of $Q_{\Lambda_1}$, and define the mapping
$$
F(\psi)= \sigma - \frac{1}{s}\int_0^s \frac{c(r)a'(r)}{r^\alpha\psi^\alpha(r)}\, {\rm d}r
$$
on $X$. The first step consists on proving that this mapping is a contraction  on $X$ equipped with the norm $\|\psi\|_X=\sup\limits_{0\leq s\leq \delta} \|\psi(s)\|$. It is easy to see that
$$
\|F(\psi)-F(\phi)\|_X\leq \Lambda_1\sup\limits_{0\leq s\leq \delta} \left\|\frac{1}{\psi^\alpha(s)} - \frac{1}{\phi^\alpha(s)}\right\|.
$$
Applying the Mean Value Theorem and the lower bound $\psi\ge\nu$, it is easy to check that
$$
\|F(\psi)-F(\phi)\|_X\leq \frac{\alpha\Lambda_1}{\nu^{1+\alpha}}\|\psi - \phi\|_X.
$$
Hence $F$  is a contraction, since
$$
\frac{\alpha\Lambda_1}{\nu^{1+\alpha}} = \frac{\alpha(\sigma-\nu)}{\nu}<1\Leftrightarrow \nu > \frac{\alpha\sigma}{\alpha+1}=z_{min}\text{ of } Q_{\Lambda_1}.
$$
Therefore, $F$ has a unique fixed point $\psi$ on $X$. Setting $\theta(s)=s\psi(s)$ provides then the existence of a unique solution of~\eqref{eq:integral_reaction_diffusion} in the class of functions satisfying $\nu s\leq \theta(s)\leq \sigma s$ for $s\in[0,\delta]$.

On the other hand, by Lemma~\ref{lemma:existence_reaction_difusion_c_positive} (i) and the extendibility result Theorem~\ref{Thm:extendibility}, any solution of~\eqref{eq:integral_reaction_diffusion} satisfies $\theta(s)\leq \sigma s$. Now, if a solution satisfies $\theta(s)\geq \gamma s$ for some $\gamma\in (\gamma^*,\nu)$, then substituting this in~\eqref{eq:integral_reaction_diffusion} one can arrive to
$$
\theta(s)\geq (\sigma - \Lambda_1/\gamma^\alpha)s\quad\text{for all }s\in(0,\delta).
$$
Iterating this process we see that $\theta(s)$ is greater than the limit of a continuous fraction of the form
$$
\sigma - \frac{\Lambda_1}{\left( \sigma - \frac{\Lambda_1}{\left( \cdots\right)^\alpha} \right)^\alpha}.
$$
Studying this limit one sees that the sequence of partial values is monotone increasing and so the limit is precisely $\nu$. Recalling again Theorem~\ref{Thm:extendibility} this implies that $\theta$ also satisfies $\theta(s)\geq \nu s$. Therefore any solution that satisfies our hypothesis for some $\gamma>\gamma^*$ must be in the class of functions were the equation is uniquely solvable.

\noindent (ii) Lemma~\ref{lemma:existence_reaction_difusion_c_positive} (i) provides the bound for $\sigma$. To prove the behaviour as $s\downarrow 0$ we  may suppose $\sigma>\sigma_A$ without loss of generality. Let $c_2(r)a'_2(r):=As^\alpha$ and apply Theorem~\ref{thm:special_case} (i) to see that for every $\rho$ equation~\eqref{eq:integral_reaction_diffusion} with $c_2a'_2$ admits a unique solution $\theta_\rho$ such that
$$
\theta_\rho (s)=z_1 s + \rho(A^{1/(\alpha + 1)} s )^{\gamma + 1} + O(s^{2\gamma +1})\quad\text{as } s\downarrow 0
$$
for a certain $\gamma>0$. It is worth mentioning that the fractional polynomial appearing here and to which we should apply Theorem~\ref{thm:special_case} is $P^+_{\sigma A^{-1/(\alpha+1)}}$, with root say $w_1$. But $w_1$ is a root of this polynomial if and only if $z_1:=A^{1/(\alpha+1)}w_1$ is a root of $Q_A$. To see this multiply the first polynomial by $A$ and apply the change of variables $w=A^{1/(\alpha+1)}z$.

Now, if $[0,\delta_\rho)$ denotes the maximal interval of existence of $\theta_\rho$ inside $[0,\delta)$, then $\theta_\rho$ is positive on $(0,\delta_\rho)$, $\delta_\rho$ depends continuously and monotonically on $\rho$ and $\delta_\rho \to 0$ as $\rho \to -\infty$. Using now the previous lemma, either
$$
\theta >\theta_\rho\quad\text{in }(0,\delta^*)\quad\text{ for some }\delta^*\in(0,\delta_\rho),
$$
or $\delta_\rho =\delta$ and $\theta\leq \theta_\rho$ in $(0,\delta)$. Since $\delta_\rho \to 0$ as $\rho \to -\infty$ we can choose a negative $\rho$ of sufficient magnitude such that the first option holds. This gives the desired result.

\noindent (iii) This part uses the same tools and ideas of the previous one, using also the maximal solution presented in Theorem~\ref{thm:special_case}. In the previous section the reader may also see how the case $\alpha=1$ translates to our general case $\alpha>0$ so we recall to \cite{Gilding-Kersner-2004} for the details on this part.

\noindent (iv) The first assertion comes from part (ii) of this lemma. For the second, choose $\gamma<z_1/A^{1/\alpha}$ and observe that $\theta_2(s):= \gamma (ca')^{1/\alpha}(s)$ is a solution of~\eqref{eq:complete_integral_equation} on $[0,\delta)$ with $\sigma_2s + b_2(s):= \theta_2(s) + s\gamma^\alpha$. Furthermore, $s\mapsto \sigma_2s + b_2(s) - \sigma s$ is nondecreasing on $[0,\delta)$, since
$$
(\sigma_2s + b_2(s) - \sigma s)' = \gamma ((ca')^{1/\alpha}(s))' + \frac{1}{\gamma^\alpha} - \sigma \geq \gamma A^{1/\alpha} + \frac{1}{\gamma^\alpha} - \sigma
$$
and this quantity is greater or equal than 0 if and only if $Q_A(z)\geq 0$ once we define $z:=\gamma A^{1/\alpha}$. Thus, we take $z< z_1$ which is equivalent to our condition over $\gamma$. Hence, by Lemma~\ref{lemma:comparison_TW_reaction_diffusion}, either $\theta>\theta_2$ in $(0,\delta^*)$ for some $\delta^*\in(0,\delta)$, or $\theta\leq \theta_2$ in $(0,\delta)$. But substituting the last option in the right-hand side of~\eqref{eq:integral_reaction_diffusion} we obtain $\theta(s)\leq (\sigma - 1/\gamma^\alpha)s<z_1 s$ for all $s\in(0,\delta)$, which contradicts part (ii) of this lemma. Therefore the first option holds and thus
$$
\theta(s)>\gamma (ca')^{1/\alpha}(s)\text{ for all }\gamma < \frac{z_1}{A^{1/\alpha}} \text{ and }  s\in(0,\delta^*)\text{ for some }\delta^*\in(0,\delta).
$$
From here the conclusion follows.

\noindent (v) The function $\theta_1(s):=\gamma((ca')(s))^{1/\alpha}$ with $\gamma:=\beta_1/B^{1/\alpha}$ is a solution of~\eqref{eq:complete_integral_equation} on $[0,\delta)$ with $\sigma_1s + b_1(s):=\theta_1(s) + s/\gamma^\alpha$. Furthermore, by ideas similar to the ones in part (iv), $s\mapsto \sigma s - \sigma_1s-b_1(s)$ is nondecreasing on $[0,\delta)$. By Lemma~\ref{lemma:comparison_TW_reaction_diffusion} either $\theta<\theta_1$ in $(0,\delta^*)$ for some $\delta^*\in(0,\delta)$, or $\theta\geq \theta_1$ in $(0,\delta)$. If this last option is true, then, substituting in~\eqref{eq:integral_reaction_diffusion}, we see that
$$
\theta(s) \geq\sigma s - \frac{s}{\gamma^\alpha}=\beta_1 \quad\text{for all }s\in(0,\delta).
$$
Applying part (i) of this lemma and Lemma~\ref{lemma:existence_reaction_difusion_c_positive} (iv) we obtain that  this solution is unique and must therefore be the maximal solution. Since by hypothesis $\theta$ is not the maximal solution the first option must hold, and the conclusion follows.
\end{proof}

The next theorem is a consequence of this lemma and Theorem~\ref{thm:bounded_support_in_general} when it is translated to the case of general semi-wavefront profiles.
\begin{Theorem}
Suppose that the conditions of Theorem~\ref{thm:existence_positive_reaction} hold. Fix $\delta\in(0,\ell)$.
\begin{itemize}
\item[\rm (i)] If
$$
\int_0^\delta \frac{a'(s)}{s^\alpha}\, {\rm d}s =\infty,
$$
then every semi-wavefront profiles decreasing to 0 is positive everywhere in its domain of definition.

\item[\rm (ii)] If
$$
\int_0^\delta \frac{a'(s)}{s^\alpha}\, {\rm d}s <\infty,
$$
then for every wave speed $\sigma>\sigma_{\lambda_1}$ there is a semi-wavefront profile decreasing to~0 whose support is bounded above. Moreover, if $ca'$ is absolutely continuous on~$[0,\delta)$, $(ca')(0)=0$, $\lim \esssup\limits_{s\to 0} ((ca')^{1/\alpha})'(s)\leq \lambda_1$ and
$$
\int_0^\delta \frac{1}{c(s)}\, {\rm d}s =\infty,
$$
then for every wave speed $\sigma>\sigma_{\lambda_1}$ there is exactly one distinct solution of this type whose support is bounded above and all other solutions of this type are positive everywhere in their domain of definition. On the contrary, if $\lim \essinf\limits_{s\to 0} (ca')(s)/s>0$ or if $ca'$ is absolutely continuous on $[0,\delta)$, $(ca')(0)=0$, $\lim \essinf\limits_{s\to 0} ((ca')^{1/\alpha})'(s)>-\infty$ and
$$
\int_0^\delta \frac{1}{c(s)}\, {\rm d}s <\infty,
$$
then for every wave speed $\sigma>\sigma_{\lambda_1}$ every solution of this type has its support bounded above.
\end{itemize}
\end{Theorem}
\subsection{Wavefronts}
\label{sect-wavefronts_reaction_diffusion}
We begin by presenting a connection between the speed $\sigma$ and the integral
\begin{equation}\label{eq:kappa}
\kappa := \int_0^\ell c(s)a'(s)\, {\rm d}s,
\end{equation}
a result that corresponds  in the semi-linear case $p=2$, $a(u)=u$ to \cite[Theorem 10.1]{Gilding-Kersner-2004}.
\begin{Theorem}\label{thm:existence_wavefronts_sigma_kappa}
Let $\ell<\infty$.
Equation~\eqref{eq:differential_reaction_diffusion} has a wavefront profile from 0 to $\ell$ with speed $\sigma$ only if one of the following holds:
\begin{itemize}
\item[\rm (a)] $\kappa >0$, $\sigma >0$ and $\displaystyle\int_u^\ell c(s)a'(s)\, {\rm d}s >0$  for all $u\in(0,\ell)$.
\item[\rm (b)] $\kappa =0$, $\sigma =0$, $\displaystyle\int_0^u c(s)a'(s)\, {\rm d}s \leq 0$ for all $u\in(0,\ell)$, and
$\displaystyle \theta(s)=G(s)$ satisfies the integrability condition on $ I$.
\item[\rm (c)] $\kappa<0$, $\sigma <0$ and $\displaystyle\int_0^u c(s)a'(s)\, {\rm d}s <0$ for all $u\in(0,\ell)$.
\end{itemize}
The necessary condition {\rm(b)} is also a sufficient condition for existence.
\end{Theorem}

\subsubsection{Fixed sign}\label{sect-fixed_sign}

Now we focus our attention on the case in which $c>0$ in $(0,\ell)$. Due to Theorem~\ref{thm:correspondance_u_l-u} our results will apply, with the needed changes, to the case in which $c<0$ in $(0,\ell)$.

\begin{Theorem}
Suppose that $\ell<\infty$ and $c>0$ in $(0,\ell)$. Set
$$
\lambda_1:=\limsup\limits_{s\downarrow 0} \left\{ \frac{1}{s}\int_0^s\frac{c(r)a'(r)}{r^\alpha} \, {\rm d}r \right\}.
$$
\begin{itemize}
\item[\rm (i)] If $\lambda_1=\infty$, then equation~\eqref{eq:differential_reaction_diffusion} has no wavefront profile from $\ell$ to 0.
\item[\rm (ii)] If $\lambda_1<\infty$, there exists a value $\sigma^*>0$ such that~\eqref{eq:differential_reaction_diffusion} has exactly one wavefront profile from $\ell$ to 0 for every wave speed $\sigma\geq\sigma^*$ and no such profile for $\sigma < \sigma^*$.
\end{itemize}
\end{Theorem}
The proof is similar to the one for  \cite[Theorem 10.5]{Gilding-Kersner-2004}

Our next step is to say something about the supports of the solutions, but prior to this we need three lemmata similar to~\cite[Lemmata 10.9, 10.19 and 10.20]{Gilding-Kersner-2004}, with again analogous proofs.
\begin{Lemma}
Let $\ell<\infty$, $c>0$ in $(0,\ell)$, $\lambda_1<\infty$ and $\sigma^*>0$ be the critical wave speed for which equation~\eqref{eq:differential_reaction_diffusion} has exactly one distinct wavefront profile from $\ell$ to 0 for every speed $\sigma\geq\sigma^*$ and no such solution for $\sigma<\sigma^*$.
\begin{itemize}
\item[\rm (i)] If $\sigma>\sigma^*$ or $\sigma\geq \sigma_{\Lambda_1}$ any solution $\theta$ of~\eqref{eq:integral_reaction_diffusion} on $\bar{I}$ satisfying $\theta(l)=0$ cannot be the maximal solution of this equation.
\item[\rm (ii)] If $\sigma=\sigma^*>\sigma_\lambda$ any solution $\theta^*$ of~\eqref{eq:integral_reaction_diffusion} on $\bar{I}$ satisfying $\theta^*(l)=0$ must be the maximal solution of this equation.
\end{itemize}
\end{Lemma}

The following is direct consequence of Lemmata~\ref{lemma:existence_reaction_difusion_c_positive} and~\ref{lemma:behaviour_near_0_reaction_diffusion_c_positive}.
\begin{Lemma}
Let $c>0$ in $(0,\delta)$ for some $\delta\in(0,\ell)$ and $\lambda_0=\lambda_1$.

\noindent{\rm (i)} If $\sigma>\sigma_\lambda$ or if $\sigma=\sigma_\lambda$ and $c(s)a'(s)/s^{\alpha}\to \lambda_1$ as $s\downarrow 0$ then the maximal solution~$\theta$ of equation~\eqref{eq:integral_reaction_diffusion} satisfies $\theta(s)/s\to z_2$, where $z_2$ is the biggest root of $Q_{\lambda_1}$.

\noindent{\rm (ii)} If $\sigma\geq\sigma_\lambda$ and $c(s)a'(s)/s^{\alpha}\to \lambda_1$ as $s\downarrow 0$ then any solution $\theta$ of equation~\eqref{eq:integral_reaction_diffusion} other than the maximal solution satisfies $\theta(s)/s\to z_1$, where $z_1$ is the smallest root of $Q_{\lambda_1}$.

\noindent{\rm (iii)} If $\sigma\geq\sigma_\lambda$, $(ca')^{1/\alpha}$ is differentiable on $[0,\delta)$, $(ca')(0)=0$ and there holds $((ca')^{1/\alpha})'\to\lambda_1^{1/\alpha}$ as $s\downarrow 0$, then any solution $\theta$ of equation~\eqref{eq:integral_reaction_diffusion} other than the maximal solution satisfies
$$
\frac{\theta(s)}{(c(s)a'(s))^{1/\alpha}}\to \frac{z_1}{\lambda_1^{1/\alpha}}.
$$
\end{Lemma}

The third one recovers some properties from negative reactions adapted to wavefronts.
\begin{Lemma}\label{lemma:behaviour_near_0_negative_kernel}
Suppose that $c\leq 0$ in $(0,\delta)$ for some $\delta\in(0,\ell)$ and that equation~\eqref{eq:integral_reaction_diffusion} has a unique solution $\theta$ on $[0,\delta)$.

\noindent{\rm (i)} If $\sigma=0$ or if $G(s)/s\to \infty$ as $s\to 0$, then $\theta(s)/G(s)\to 1$ as $s\to 0$.

\noindent{\rm (ii)} If $\sigma >0$ and $G(s)/s\to \mu$ as $s\to 0$ for some $\mu\in[0,\infty)$, then $\theta(s)/G(s)\to \mu^*$ as $s\to 0$, being $\mu^*$ a positive constant dependant on $\mu$.

\noindent{\rm (iii)} If $c<0$ in $(0,\delta)$, $ca'$ is differentiable on $(0,\delta)$, $(ca')(0)=0$, $(|ca'|^{1/\alpha})'(s)\to \mu$ as $s\downarrow 0$ for some $\mu\in [0,\infty)$, then $\theta(s)/|c(s)a'(s)|^{1/\alpha}\to \overline{\mu}$ as $s\to 0$, being $\overline{\mu}$ a positive constant dependant on $\mu$.
\end{Lemma}

With this lemmata on our hands, we next show the theorem about the support of the solutions. Its proof is similar to  \cite[Theorem 10.21]{Gilding-Kersner-2004}.
\begin{Theorem}\label{thm:suport_for_wavefronts_fixed_sign}
Suppose that $\ell<\infty$, $c>0$ in $(0,\ell)$ and $\lambda_1<\infty$. Let $\sigma^*>0$ be the critical wave speed for which equation~\eqref{eq:differential_reaction_diffusion} has exactly one distinct wavefront profile from $\ell$ to 0 for every speed $\sigma\geq\sigma^*$ and no such solution for $\sigma<\sigma^*$.
\begin{itemize}
\item[\rm (i)] Suppose in addition that $(ca')^{1/\alpha}$ is differentiable on $[0,\delta]$ for some $\delta\in(0,\ell)$, $(ca')(0)=0$ and $((ca')^{1/\alpha})'(u)\to ((ca')^{1/\alpha})'(0)$ as $u\downarrow 0$. Then the following alternatives are mutually exclusive.
\begin{itemize}
\item[\rm (a)] Every wavefront profile from $\ell$ to 0 satisfies~\eqref{eq:suport_bouded_above}. This occurs if and only if
$$
\int_0^\delta \frac{1}{c(s)}\, {\rm d}s <\infty.
$$
\item[\rm (b)] Every wavefront profile from $\ell$ to 0 with wave speed $\sigma^*$ satisfies~\eqref{eq:suport_bouded_above}, while every such solution with speed $\sigma>\sigma^*$ satisfies~\eqref{eq:suport_not_bouded_above}. This occurs if and only if
$$
\int_0^\delta \frac{a'(s)}{s^\alpha}\, {\rm d}s<\int_0^\delta \frac{1}{c(s)}\, {\rm d}s =\infty.
$$
\item[\rm (c)] Every wavefront profile from $\ell$ to 0 satisfies~\eqref{eq:suport_not_bouded_above}. This occurs if and only if
$$
\int_0^\delta \frac{a'(s)}{s^\alpha}\, {\rm d}s =\infty.
$$
\end{itemize}

\item[\rm (ii)] Suppose furthermore that $(ca')^{1/\alpha}$ is differentiable on $[\ell-\delta,\ell]$ for some $\delta\in(0,\ell)$, $(ca')(\ell)=0$ and $((ca')^{1/\alpha})'(u)\to ((ca')^{1/\alpha})'(\ell)$ as $u\uparrow \ell$. Then the following alternatives are mutually exclusive.
\begin{itemize}
\item[\rm (a)] Every wavefront profile from $\ell$ to 0 satisfies~\eqref{eq:suport_bouded_below}. This occurs if and only if
$$
\int_{\ell-\delta}^\ell \frac{1}{c(s)}\, {\rm d}s <\infty.
$$
\item[\rm (b)] Every wavefront profile from $\ell$ to 0 satisfies~\eqref{eq:suport_not_bouded_below}. This occurs if and only if
$$
\int_{\ell-\delta}^\ell \frac{1}{c(s)}\, {\rm d}s =\infty.
$$
\end{itemize}
\end{itemize}
\end{Theorem}

We close this section studying the behaviour of our wavefront profiles when they approach the critical values 0 and $\ell$. This result is analogous to \cite[Theorem 10.22]{Gilding-Kersner-2004} though we find important to sketch the behaviour of the parameter $\lambda_1$ in this case. When $\alpha=1$ the work of Gilding and Kersner shows how $\lambda_1$ becomes $(ca')'(0)$ but our case is not so tidy. In fact tracing back the definition of $\lambda_1$ we see that in our case, whenever $\lambda_1<\infty$, we have that this value becomes $[((ca')^{1/\alpha})'(0)]^{\alpha}$.
\begin{Theorem}\label{thm:behaviour_near_0_wavefront_fixed_sign}
Let $V$ denote a wavefront profile of equation~\eqref{eq:differential_reaction_diffusion} from $\ell$ to 0 with wave speed $\sigma$.
\begin{itemize}
\item[\rm (i)] Suppose that the conditions of Theorem~\ref{thm:suport_for_wavefronts_fixed_sign} {\rm (i)} hold and define
\begin{equation}
\label{eq:def.xi^*}
\xi^*:=\sup\{\xi\in\mathbb{R}: V(\xi)>0\}.
\end{equation}
If $\sigma=\sigma^*$, then $\frac{(a(V))'(\xi)}{V^\alpha(\xi)} \to -z_2^\alpha$ as $\xi\uparrow\xi^*$,
where $z_2$ is the biggest root of $Q_\gamma$, with $\gamma=[((ca')^{1/\alpha})'(0)]^{\alpha}$,
whereas if $\sigma>\sigma^*$, then $\frac{V'(\xi)}{c(V(\xi))} \to -\frac{z_1^\alpha}{((ca')^{1/\alpha})'(0)}$ as $\xi\uparrow\xi^*$,
where $z_1$ is the smallest root of the \textit{fractional polynomial} $Q_\gamma$.

\item[\rm (ii)] Suppose that the conditions of Theorem~\ref{thm:suport_for_wavefronts_fixed_sign} {\rm (ii)} hold and define
\begin{equation}
\label{eq:def.xi_*}
\xi_*:=\inf\{\xi\in\mathbb{R}: V(\xi)<\ell\}.
\end{equation}
Then, $\frac{V'(\xi)}{c(V(\xi))} \to -\frac{z_0^\alpha}{((ca')^{1/\alpha})'(\ell)}$ as $\xi\downarrow\xi_*$,
where $z_0$ is the smallest root of $Q_\gamma$ with $\gamma=[((ca')^{1/\alpha})'(\ell)]^{\alpha}$.
\end{itemize}
\end{Theorem}

\subsubsection{One sign change}\label{sect-one_sign_change}

It is time to close up this work and we do so by studying the properties of wavefront profiles appearing when the reaction term has one sign change, in essence, $c\leq 0$ in  $[0,a]$ and $c\geq 0$ in $[a,\ell]$, for some $a\in (0,\ell)$. In the opposite case, when the reaction starts being non-negative and ends being non-positive,  Theorem~\ref{thm:existence_wavefronts_sigma_kappa} prevents the existence of wavefront profiles connecting $\ell$ with 0. Therefore, we stick to the first case, and to work with it we define the value
$$
H(s):=\left|(\alpha+1) \int_s^\ell c(r)a'(r)\, {\rm d}r \right|^{1/(\alpha+1)},
$$
which plays a similar role as the function $G$ but coming from the value $\ell$. We also recall the equation
\begin{equation}\label{eq:wavefront_Theta_mayuscula}
\Theta(s)=-\sigma s + \tilde{b}(s) - \int_0^s\frac{\tilde{c}(r)\tilde{a}'(r)}{\Theta^\alpha(r)} \, {\rm d}r
\end{equation}
from Theorem~\ref{thm:correspondance_theta_Theta}, and the reader may have already guessed that we are going to use it to split the analysis of our equation in two halves separated by the value $u=a$. This creates a new difficulty, because during the study of the existence of wavefront profiles we would like to paste the two halves together at some point and with the same speed of propagation. For the sake of brevity, we will omit the two lemmata that provide these in our work and just cite the correspondent ones in \cite{Gilding-Kersner-2004}, since they are analogous. Continuity of wavefront profiles with respect to the parameter $\sigma$ is stated in \cite[Lemma 8.5]{Gilding-Kersner-2004}, and the possibility of sticking both halves back together in \cite[Lemma 8.9]{Gilding-Kersner-2004}.

\begin{Theorem}\label{thm:existence_wavefront_one_sign_change}
Suppose that $c\leq 0$ in $[0,a]$ and $c\geq 0$ in $[a,\ell]$, for some $a\in(0,\ell)$, $\ell<\infty$. Let $\kappa$ be defined as in~\eqref{eq:kappa}. Suppose also that one of the following holds:
\begin{itemize}
\item[\rm (i)] $\kappa>0$ and $c> 0$ in  $(a,\ell)$;

\item[\rm (ii)] $\kappa=0$, $G(u)> 0$ in $(0,a)$ and $H(u)> 0$ in $(a,\ell)$;

\item[\rm (iii)] $\kappa<0$ and $c< 0$ in $(0,a)$.
\end{itemize}
Then there exists a real number $\sigma^*$ for which equation~\eqref{eq:differential_reaction_diffusion} has exactly one distinct wavefront profile from $\ell$ to 0 with speed $\sigma^*$ and no such solution for any other wave speed.
\end{Theorem}
\begin{proof}
Recalling Theorem~\ref{thm:formal_uniqueness}, let $S_0$ denote the set of values $\sigma$ for which equation~\eqref{eq:integral_reaction_diffusion} has a unique solution $\theta(s;\sigma)$ on $[0,a]$, $S_1$ denote the set of values for which equation~\eqref{eq:wavefront_Theta_mayuscula} has a unique solution $\Theta(s;\sigma)$ on $[0,\ell-a]$  and $S=S_0\cap S_1$. Considering $\sigma=0$ in both equations it can be seen that $G(s)$ solves~\eqref{eq:integral_reaction_diffusion} on $[0,a]$ and $H(\ell-s)$ solves~\eqref{eq:wavefront_Theta_mayuscula} on $[0,l-a]$.  By the analogous of \cite[Lemma 8.5]{Gilding-Kersner-2004}, $S_0$ is an interval containing $[0,\infty)$ with $\theta(s;0)=G(s)$ and $S_1$ is an interval containing $(-\infty,0]$ with $\Theta(s;0)=H(\ell-s)$. Furthermore, $\theta(a;0)\to \infty$ as $\sigma\to\infty$ and $\Theta(a;0)\to \infty$ as $\sigma\to-\infty$.

Let us define, for $\sigma\in S$, $F(\sigma):= (\theta(a;\sigma))^{(\alpha+1)/2}-(\Theta(l-a;\sigma))^{(\alpha+1)/2}$. Again by  \cite[Lemma 8.5]{Gilding-Kersner-2004}, $F$ is a continuous function of $\sigma\in S$ and
$$
F(0)= (G(a))^{(\alpha+1)/2}-(H(a))^{(\alpha+1)/2} = -\frac{\kappa}{(G(a))^{(\alpha+1)/2}+(H(a))^{(\alpha+1)/2}}.
$$
Now by Theorem~\ref{thm:correspondencia_ecuaciones_integral_y_diferencial} and \cite[Lemma 8.9]{Gilding-Kersner-2004}, to prove the existence of a speed for which our equation admits a wavefront profile from $\ell$ to 0 it suffices to show that there is a value $\sigma^*\in S$ such that $F(\sigma^*)=0$, $\theta(\cdot;\sigma^*)>0$ on $(0,a]$ and $\Theta(\cdot;\sigma^*)>0$ on $(0,\ell-a]$. In order to do so we distinguish the three cases that we stated in the theorem.

\noindent (i) In this case, by Theorem~\ref{thm:uniqueness_positive_kernel_or_0} applied to equation~\eqref{eq:wavefront_Theta_mayuscula}, $S_1=(-\infty,\infty)$ and therefore $S$ contains $[0,\infty)$. Furthermore, $F(0)<0$ and $F(\sigma)\to \infty$ as $\sigma\to\infty$, so there must exist a value $\sigma^*>0$ such that $F(\sigma^*)=0$. Moreover, $\theta(\cdot;\sigma^*)\geq \sigma^* s$ for all $s\in(0,a]$ an again by Theorem~\ref{thm:uniqueness_positive_kernel_or_0} applied to equation~\eqref{eq:wavefront_Theta_mayuscula}, $\Theta(\cdot;\sigma^*)>0$ on $(0,\ell-a]$.

\noindent (ii) By hypothesis, it is enough to take $\sigma^*=0$.

\noindent (iii) This case is analogous to (i), but now $S\supset(-\infty,0]$, $F(0)>0$ and $F(\sigma)\to -\infty$ as $\sigma\to-\infty$.

Let us see that the speed $\sigma^*$ is unique. Since $\theta(a;\sigma^*)=\Theta(\ell-a;\sigma^*)>0$, by Theorem~\ref{thm:formal_comparison_maximal} applied to equations~\eqref{eq:integral_reaction_diffusion} and~\eqref{eq:wavefront_Theta_mayuscula} there holds $F(\sigma)<0$ for all $\sigma\in S\cap (-\infty,\sigma^*)$ and $F(\sigma)>0$ for all $\sigma\in S\cap (\sigma^*, \infty)$, so again by Theorem~\ref{thm:correspondencia_ecuaciones_integral_y_diferencial} and \cite[Lemma 8.9]{Gilding-Kersner-2004} the value $\sigma^*$ is the unique speed for which our wavefront profile exists. The distinctiveness of the wavefront profile follows from Corollary~\ref{thm:es_el_2.25_del_libro} when adapted to wavefronts.
\end{proof}

We finish by studying the behaviour close to the critical values $S=0,\ell$ and the support of the solutions, as we did in the previous section.
\begin{Theorem}\label{thm:suport_for_wavefronts_one_sign_change}
Suppose that $\ell<\infty$ and let $V$ denote the wavefront profile from $\ell$ to 0 of equation~\eqref{eq:differential_reaction_diffusion} with wave speed $\sigma$.
\begin{itemize}
\item[\rm (i)] Suppose also that $c\leq 0$ on $(0,\delta]$ for some $\delta\in (0,\ell/2)$, $(ca')^{1/\alpha}$ is differentiable on $[0,\delta]$, $(ca')(0)=0$ and $((ca')^{1/\alpha})'(u)\to ((ca')^{1/\alpha})'(0)$ as $u\downarrow 0$.
\begin{itemize}
\item[\rm (a)] If $\kappa>0$ then $V$ satisfies~\eqref{eq:suport_bouded_above} if and only if
$$
\int_0^\delta \frac{a'(s)}{s^\alpha}\, {\rm d}s <\infty.
$$
\item[\rm (b)] If $\kappa=0$ and $G>0$ on $(0,\delta]$ then $V$ satisfies~\eqref{eq:suport_bouded_above} if and only if
$$
\int_0^\delta \frac{a'(s)}{G^\alpha(s)}\, {\rm d}s<\infty.
$$
\item[\rm (c)] If $\kappa<0$ and $c<0$ on $(0,\delta]$ then $V$ satisfies~\eqref{eq:suport_bouded_above} if and only if
$$
\int_0^\delta \frac{1}{|c(s)|}\, {\rm d}s<\infty.
$$
\end{itemize}

\item[\rm (ii)] Suppose that $c\geq 0$ on $[\ell-\delta,\ell)$ for some $\delta\in(0,\ell/2)$, $(ca')^{1/\alpha}$ is differentiable on $[\ell-\delta, \ell]$, $(ca')(\ell)=0$ and $((ca')^{1/\alpha})'(u)\to ((ca')^{1/\alpha})'(\ell)$ as $u\uparrow \ell$.
\begin{itemize}
\item[\rm (a)] If $\kappa>0$ and $c>0$ on $[\ell-\delta,\ell)$ then $V$ satisfies~\eqref{eq:suport_bouded_below} if and only if
$$
\int_{\ell-\delta}^\ell\frac{1}{c(s)}\, {\rm d}s <\infty.
$$
\item[\rm (b)] If $\kappa=0$ and $H>0$ on $[\ell-\delta,\ell)$ then $V$ satisfies~\eqref{eq:suport_bouded_below} if and only if
$$
\int_{\ell-\delta}^\ell \frac{a'(s)}{H^\alpha(s)}\, {\rm d}s<\infty.
$$
\item[\rm (c)] If $\kappa<0$ then $V$ satisfies~\eqref{eq:suport_bouded_above} if and only if
$$
\int_{\ell-\delta}^\ell \frac{a'(s)}{(\ell-s)^\alpha}\, {\rm d}s<\infty.
$$
\end{itemize}
\end{itemize}
\end{Theorem}

This theorem is a corollary of the following, which is at the same time consequence of Theorems~\ref{thm:es_el_2.25_del_libro},  \ref{thm:existence_wavefronts_sigma_kappa}, \ref{thm:existence_wavefront_one_sign_change}, \ref{thm:correspondance_theta_Theta} and Lemma~\ref{lemma:behaviour_near_0_negative_kernel}.
\begin{Theorem}
Let $V$ denote a wavefront profile of equation~\eqref{eq:differential_reaction_diffusion} from $\ell$ to 0 with wave speed $\sigma$. Define $\xi^*$ and $\xi_*$ as in~\eqref{eq:def.xi^*}--\eqref{eq:def.xi_*}.

\noindent{\rm (i)} Suppose that the conditions of Theorem~\ref{thm:suport_for_wavefronts_one_sign_change} {\rm (i)} hold.
\begin{itemize}
\item[\rm (a)] If $\kappa>0$ then
$
\frac{(a(V))'(\xi)}{V^\alpha(\xi)} \to -z_2^\alpha$ as $\xi\uparrow\xi^*$,
where $z_2$ is the biggest root of $Q_\gamma$ with $\gamma=[((ca')^{1/\alpha})'(0)]^{\alpha}$.

\item[\rm (b)] If $\kappa=0$ and $G>0$ on $(0,\delta]$ then
$
\frac{(a(V))'(\xi)}{G(V^\alpha(\xi))} \to -1$ as $\xi\uparrow\xi^*$.

\item[\rm (c)] If $\kappa<0$ and $c<0$ on $(0,\delta]$ then
$
\frac{V'(\xi)}{|c(V(\xi))|} \to \frac{z_2^\alpha}{((ca')^{1/\alpha})'(0)}\quad\text{as }\xi\uparrow\xi^*$.
\end{itemize}

\noindent{\rm (ii)} Suppose that the conditions of Theorem~\ref{thm:suport_for_wavefronts_fixed_sign} {\rm (ii)} hold.
\begin{itemize}
\item[\rm (a)] If $\kappa>0$ and $c>0$ on $[\ell-\delta,\ell)$ then
$
\frac{V'(\xi)}{|c(V(\xi))|} \to -\frac{z_1^\alpha}{((ca')^{1/\alpha})'(0)}$ as $\xi\downarrow\xi_*$,
where $z_1$ is the smallest root of $Q_\gamma$ with $\gamma=[((ca')^{1/\alpha})'(\ell)]^{\alpha}$.

\item[\rm (b)] If $\kappa=0$ and $H>0$ on $[\ell-\delta,\ell)$ then
$\frac{(a(V))'(\xi)}{H(V^\alpha(\xi))} \to -1$ as $\xi\downarrow\xi_*$.
\item[\rm (c)] If $\kappa<0$ then $\frac{(a(V))'(\xi)}{(\ell-V(\xi))^\alpha} \to z_1^\alpha$ as $\xi\downarrow\xi_*$.
\end{itemize}
\end{Theorem}


\

\noindent{\large \textbf{Acknowledgments}}

\noindent  The author wants to thank his thesis advisor, Fernando Quir\'os, for his help and encouragement during the realization of this work.



\

\noindent\textbf{Address:}

\noindent\textsc{A. G\'arriz: } Departamento de Matem\'{a}ticas, Universidad
Aut\'{o}noma de Madrid, 28049 Madrid, Spain.  (e-mail: alejandro.garriz@estudiante.uam.es).

\end{document}